\documentclass[11pt]{article}
\usepackage{amssymb}
\usepackage{amsmath}
\usepackage{booktabs}
\usepackage{graphicx}
\usepackage{theorem}

\newtheorem{thm}{Theorem}
\newtheorem{lem}{Lemma}
\newtheorem{cor}{Corollary}
\newtheorem{prp}{Proposition}
\theorembodyfont{\rmfamily}

\allowdisplaybreaks[3]
\makeatletter
\newcommand{\Proofname}{Proof}
\makeatother
\makeatletter
\def\BOXSYMBOL{\RIfM@\bgroup\else$\bgroup\aftergroup$\fi
  \vcenter{\hrule\hbox{\vrule height.85em\kern.6em\vrule}\hrule}\egroup}
\makeatother
\newcommand{\BOX}{%
  \ifmmode\else\leavevmode\unskip\penalty9999\hbox{}\nobreak\hfill\fi
  \quad\hbox{\BOXSYMBOL}}
\def\al{{\alpha}}
\def\be{{\beta}}

\def\la{{\lambda}}

\def\fai{{\varphi}}

\def\bbe{{\text{\boldmath $\beta$}}}

\def\bep{{\text{\boldmath $\varepsilon$}}}

\def\bth{{\text{\boldmath $\theta$}}}

\def\bpsi{{\text{\boldmath $\psi$}}}

\def\bbeh{{\widehat \bbe}}

\def\bthh{{\widehat \bth}}

\def\btht{{\widetilde \bth}}

\def\Ga{{\Gamma}}

\def\bPsi{{\text{\boldmath $\Psi$}}}

\def\bPsih{{\widehat \bPsi}}

\def\g{{\text{\boldmath $g$}}}

\def\u{{\text{\boldmath $u$}}}
\def\v{{\text{\boldmath $v$}}}

\def\y{{\text{\boldmath $y$}}}
\def\z{{\text{\boldmath $z$}}}
\def\A{{\text{\boldmath $A$}}}
\def\B{{\text{\boldmath $B$}}}
\def\C{{\text{\boldmath $C$}}}
\def\D{{\text{\boldmath $D$}}}

\def\G{{\text{\boldmath $G$}}}
\def\H{{\text{\boldmath $H$}}}
\def\I{{\text{\boldmath $I$}}}

\def\K{{\text{\boldmath $K$}}}

\def\P{{\text{\boldmath $P$}}}

\def\T{{\text{\boldmath $T$}}}
\def\U{{\text{\boldmath $U$}}}
\def\V{{\text{\boldmath $V$}}}

\def\X{{\text{\boldmath $X$}}}

\def\ah{{\hat a}}
\def\bh{{\hat b}}

\def\Nc{{\cal N}}

\def\Re{{\mathbb{R}}}
\def\tr{{\rm tr\,}}
\def\diag{{\rm diag\,}}

\def\[{{\text{\boldmath $[$}}}
\def\]{{\text{\boldmath $]$}}}
\def\et{{\it et\, al.}}
\def\zero{{\bf\text{\boldmath $0$}}}

\def\|{{\,|\,}}

\def\/{{\Bigr/\!\!}}

\def\1r{{\rm (1)}}
\def\2r{{\rm (2)}}
\def\3r{{\rm (3)}}
\def\4r{{\rm (4)}}
\def\5r{{\rm (5)}}

\def\non{{\nonumber}}

\def\zero{{\bf\text{\boldmath $0$}}}

\begin{document}
\title{Corrected Empirical Bayes Confidence Region in a Multivariate Fay-Herriot Model}

\author{
Tsubasa Ito\footnote{Graduate School of Economics, University of Tokyo, 7-3-1 Hongo, Bunkyo-ku, Tokyo 113-0033, JAPAN. {E-Mail: tsubasa$\_$ito.0710@gmail.com}}
and
Tatsuya Kubokawa\footnote{Faculty of Economics, University of Tokyo, 7-3-1 Hongo, Bunkyo-ku, Tokyo 113-0033, JAPAN. \newline{E-Mail: tatsuya@e.u-tokyo.ac.jp }} 
}
\date{}
\maketitle
\begin{abstract}
In the small area estimation, the empirical best linear unbiased predictor (EBLUP) in the linear mixed model is useful because it gives a stable estimate for a mean of  a small area.
For measuring uncertainty of EBLUP, much of research is focused on second-order unbiased estimation of mean squared prediction errors in the univariate case.
In this paper, we consider the multivariate Fay-Herriot model where  the covariance matrix of random effects is fully unknown, and obtain a confidence reagion of the small area mean that is based on the Mahalanobis distance centered around EBLUP and is second order correct.
A positive-definite, consistent and second-order unbiased estimator of the covariance matrix of the random effects is also suggested.
The performance is investigated through simulation study.

\par\vspace{4mm}

\end{abstract}

\section{Introduction}

Mixed effects models and their model-based estimators have been recognized as useful methods in statistical inference.
In particular, small area estimation is an important application of mixed effects models.
Although direct design-based estimates for small area means have large standard errors because of small sample sizes from small areas, the empirical best linear unbiased predictors (EBLUP) induced from mixed effects models provide reliable estimates by ^^ ^^ borrowing strength" from neighboring areas and by using data of auxiliary variables.
Such a model-based method for small area estimation has been studied extensively and actively from both theoretical and applied aspects, mostly for handling univariate survey data.
For comprehensive reviews of small area estimation, see  Ghosh and Rao (1994), Datta and Ghosh (2012), Pfeffermann (2013) and Rao and Molina (2015).

\medskip
When multivariate data with correlations are observed from small areas for estimating multi-dimensional characteristics, like poverty and unemployment indicators, Fay (1987) suggested a multivariate extension of the univariate Fay-Herriot model, called a multivariate Fay-Herriot model, to produce reliable estimates of median incomes for four-, three- and five-person families.
Fuller and Harter (1987) also considered a multivariate modeling for estimating a finite population mean vector.
Datta, Day and Basawa (1999) provided unified theories in empirical linear unbiased prediction or empirical Bayes estimation in general multivariate mixed linear models. 
Datta, Day and Maiti (1998) suggested a hierarchical Bayesian approach to multivariate small area estimation.
Datta, $\et$ (1999) showed the interesting result that the multivariate modeling produces more efficient predictors than the conventional univariate modeling.
Porter, Wikle and Holan (2015) used the multivariate Fay-Herriot model for modeling spatial data.
Ngaruye, von Rosen and Singull (2016) applied a  multivariate mixed linear model to crop yield estimation in Rwanda.

\medskip
Although Datta, $\et$ (1999) developed the general and unified theories concerning the empirical best linear unbiased predictors (EBLUP) and their uncertainty,  it is definitely more helpful and useful to provide concrete forms with closed expressions for EBLUP, the second-order approximation of the mean squared error matrix (MSEM) and the second-order unbiased estimator of the mean squared error matrix.
Recently, Benavent and Morales (2016) treated the multivariate Fay-Herriot model with the covariance matrix of random effects depending on unknown parameters.
As a structure in the covariance matrix, they considered diagonal, AR(1) and the related structures and employed the residual maximum likelihood (REML) method for estimating the unknown parameters embedded in the covariance matrix.
A second-order approximation and estimation of the MESM were also derived.
However, they did not concern about the construction of confidence regions for small area means.

\medskip
Confidence regions are more useful for measuring uncertainty of EBLUP, but there is no literature about confidence regions for multivariate small area estimation problems to the best of our knowledge.
Naive confidence regions can be constructed easily by using the Bayes estimators of small area means and their MSEM.
As is the case in the univariate small area estimation problem, the coverage probability of the naive methods cannot be guaranteed to be greater than or equal to the nominal confidence coefficient $1-\al$.
Recently, in the univariate Fay-Herriot model, Diao, Smith, Datta, Maiti and Opsomer (2014) constructed closed-form confidence intervals whose coverage probability is identical to the nominal confidence coefficient up to the second-order for small area means  under the normality assumption.

\medskip
In this paper, we consider the problem of costructing confidence regions for small area mean vectors in the multivariate Fay-Herriot model where the covariance matrix of random effects is fully unknown.
Although this is a multivaliate extension of Diao et al. (2014), we are faced with two difficulties: 
One is how to construct a confidence region on the multi-dimensional space, and the other is how to construct a positive-definite and consistent estimator of the covariance matrix of random effects.
We here consider a confidence region based on the Mahalanobis distance centerd around EBLUP, and use the asymptotic expansion of the characteristic function of this distance to approximate the coverage probability based on the chi-square distributions.
We obtain the correction term in a closed form, and provide the confidence region that is second order correct.
Concerning the estimation of the covariance matrix, the Prasad-Rao type estimator with non-negative definite modification can be given in a closed form by the moment method.
When the covariance matrix is estimated with the zero matrix or a singular matrix close to the zero matrix, however, the correction term becomes instable in the confidence region.
This fact is well known in the univariate confidence interval.
Thus, we suggest a new method for obtaining a positive-definite and sencon-order unbiased estimator of the covariance matrix.
Moreover, we extend our results to construction of corrected confidence regions for the difference of two small area mean vectors.
Another approach to construction of corrected confidence regions is the bootstrap method which needs heavy burden in computation.
Because the corrected confidence region suggested here is provided in closed forms, it is easy to implement, which is a merit of our method.

\medskip
The paper is organized as follows: 
Section \ref{sec:EBLUP} introduces the multivariate Fay-Herriot model and gives the EBLUP and its prediction risk approximation.
In section \ref{sec:cr}, our proposed confidence region is derived.
Section \ref{sec:psi} gives the Prasad-Rao type estimator of the covariance matrix of the random effects and its positive-definite modification with second-order unbiasedness and consistency.
In section \ref{sec:crd}, the extension to the confidence regions for the difference of two small area means  is described.
The performances of our proposed methods are investigated in Section \ref{sec:sim}.

\section{Multivariate Fay-Herriot Model and Empirical Best Linear Unbiased Predictor}
\label{sec:EBLUP}

Suppose that area-level data $(\y_1, \X_1), \ldots, (\y_m, \X_m)$ are observed, where $m$ is the number of small areas, $\y_i$ is a $k$-variate vector of direct survey estimates and $\X_i$ is a $k\times s$ matrix of covariates associated with $\y_i$ for the $i$-th area.
The multivariate Fay-Herriot model suggested by Fay (1987) is described as
\begin{equation}
\y_i = \X_i \bbe + \v_i + \bep_i,
\quad i=1, \ldots, m,
\label{eqn:MFH}
\end{equation}
where $\bbe$ is an $s$-variate vector of unknown regression coefficients, $\v_i$ is a $k$-variate vector of random effects depending on the $i$-th area and $\bep_i$ is a $k$-variate vector of sampling errors.
It is assumed that $\v_i$ and $\bep_i$ are mutually independently distributed as
\begin{equation}
\v_i \sim \Nc_k (\zero, \bPsi)\quad \text{and}\quad \bep_i\sim\Nc_k(\zero, \D_i),
\label{eqn:MFH0}
\end{equation}
where $\bPsi$ is a $k\times k$ unknown and nonsingular covariance matrix and $\D_1, \ldots, \D_m$ are $k\times k$ known covariance matrices.
This is a multivariate extension of the so-called Fay-Herriot model suggested by Fay and Herriot (1979).
Letting $\bth_i=\X_i\bbe+\v_i$ for $i=1, \ldots, m$, we can rewrite the model given in (\ref{eqn:MFH}) and (\ref{eqn:MFH0}) as
\begin{equation}
\begin{split}
\y_i \mid \bth_i \sim& \Nc_k(\bth_i, \D_i),\\
\bth_i \sim& \Nc_k(\X_i\bbe, \bPsi),
\end{split}
\label{eqn:Bmodel}
\end{equation}
for $i=1, \ldots, m$.
Thus, the multivariate Fay-Herriot model is interpreted as the Bayes model with the prior distribution of $\bth_i$.
It may be convenient to express model (\ref{eqn:MFH}) in a matrix form.
Let $\y=(\y_1^\top, \ldots, \y_m^\top)^\top$, $\X=(\X_1^\top, \ldots, \X_m^\top)^\top$, $\v=(\v_1^\top, \ldots, \v_m^\top)^\top$ and $\bep=(\bep_1^\top, \ldots, \bep_m^\top)^\top$.
Then, model (\ref{eqn:MFH}) is expressed as
\begin{equation}
\y=\X\bbe + \v + \bep,
\label{eqn:MFH1}
\end{equation}
where $\v\sim \Nc_{km}(\zero, \I_m \otimes \bPsi)$ and $\bep\sim\Nc_{km}(\zero, \D)$ for $\D=\text{block\ diag}(\D_1, \ldots, \D_m)$.
Throughout the paper, it is assumed that $\X$ is of full rank.

\medskip
For example, we consider the crop data of Battese, Harter and Fuller (1988), who analyze the data in the nested error regression model.
For the $i$-th county, let $y_{i1}$ and $y_{i2}$ be survey data of average areas of corn and soybean, respectively.
Also let $x_{i1}$ and $x_{i2}$ be satellite data of average areas of corn and soybean, respectively.
In this case, $\y_i$, $\X_i$ and $\bbe$ correspond to
$$
\y_i=(y_{i1}, y_{i2})^\top, \quad \X_i=\begin{pmatrix} 1 & x_{i1} & x_{i2} & 0 & 0 & 0\\ 0&0&0&1 & x_{i1} & x_{i2} \end{pmatrix},\quad
\bbe=(\be_1, \ldots, \be_6)^\top
$$
for $k=2$ and $s=6$.
Battese, $\et$ (1988) applied a univariate nested error regression model for each of $y_{i1}$ and $y_{i2}$, while we can use the multivariate model (\ref{eqn:MFH}) for analyzing both data simultaneously.

\medskip
In this paper, we want to construct a confidence region of $\bth_a$ for the $a$-th area.
To this end, we begin by deriving the Bayes estimator of $\bth_a$.
The posterior distribution of $\bth_i$ given $\y_i$ and the marginal distribution of $\y_i$ are
\begin{equation}
\begin{split}
\bth_i \mid \y_i \sim & \Nc_k(\btht_a(\bbe, \bPsi), (\bPsi^{-1}+\D_i^{-1})^{-1}),\\
\y_i \sim& \Nc_k(\X_i\bbe, \bPsi+\D_i),
\end{split}
\quad i=1, \ldots, m,
\label{eqn:post}
\end{equation}
where
\begin{align*}
\btht_a(\bbe, \bPsi)=&
\X_i^\top\bbe+\bPsi(\bPsi+\D_i)^{-1}(\y_i - \X_i\bbe)\\
=&
\y_i  - \D_i(\bPsi+\D_i)^{-1}(\y_i-\X_i\bbe),
\end{align*}
which is the Bayes estimator of $\bth_i$.

\medskip
When $\bPsi$ is known, the maximum likelihood estimator or generalized least squares estimator of $\bbe$ is
\begin{align}
\bbeh(\bPsi) =& 
\{\X^\top (\I_m\otimes \bPsi + \D)^{-1}\X\}^{-1}\X^\top (\I_m\otimes \bPsi + \D)^{-1}\y\non\\
=&\Big\{ \sum_{i=1}^m \X_i^\top(\bPsi+\D_i)^{-1}\X_i\Big\}^{-1} \sum_{i=1}^m \X_i^\top(\bPsi+\D_i)^{-1}\y_i.
\label{eqn:beh}
\end{align}
Substituting $\bbeh(\bPsi)$ into $\btht_a(\bbe, \bPsi)$ yields the empirical Bayes estimator
\begin{equation}
\bthh_a(\bPsi) = \y_a  - \D_a(\bPsi+\D_a)^{-1}\big\{\y_a-\X_a\bbeh(\bPsi)\big\}.
\label{eqn:Bayes}
\end{equation}
Datta, $\et$ (1999) showed that $\bthh_a(\bPsi)$ is the best linear unbiased predictor (BLUP) of $\bth_a$.
It can be also demonstrated that $\bthh_a(\bPsi)$ is the Bayes estimator against the uniform prior distribution of $\bbe$ as well as the empirical Bayes estimator as shown above, which is called the Bayes empirical Bayes estimator.

\medskip
Because $\bPsi$ is unknown, we need to estimate the covariance matrix $\bPsi$.
Estimators used in the univariate case are the ANOVA type estimator given by Prasad and Rao (1990), the Fay-Herriot estimator suggested by Fay and Herriot (1979), and the ML and REML methods used in Datta and Lahiri (2000).
Corresponding to the univariate case, we consider the general class of estimators $\bPsih$ of $\bPsi$ which satisfy the following conditions:

\medskip
(H1)\ $\bPsih$ is an even function of $\y$ ; $\bPsih(\y)=\bPsih(-\y)$

\smallskip
(H2)\ $\bPsih$ is a translation invariant function ; $\bPsih(\y+\X\T)=\bPsih(\y)$ for any $\T\in\Re^s$ and all $\y$.

\medskip
\noindent
The modified Prasad-Rao estimator suggested later in this paper and the ML method satisfy these conditions.
We replace $\bPsi$ in $\bthh_a(\bPsi)$ with the estimator $\bPsih$, and the resulting empirical Bayes (EB) estimator is
\begin{equation}
\bthh_a^{EB} = \bthh_a(\bPsih)=\y_a  - \D_a(\bPsih+\D_a)^{-1}\big\{\y_a-\X_a\bbeh(\bPsih)\big\}.
\label{eqn:EB}
\end{equation}
This is also interpreted as the empirical best linear unbiased predictor (EBLUP) in the context of the linear mixed models.

\medskip
For evaluating the uncertainty of $\bthh_a^{EB}$, we prepare three lemmas.

\begin{lem}
\label{lem:1}
$\bbeh(\bPsi)$ is independent of $\P \y$ for $\P=\I-\X(\X^\top\X)^{-1}\X$.
Also, $\bthh_a^{EB}-\bthh_a(\bPsi)$ is a function of $\P \y$, and independent of $\bbeh(\bPsi)$.
\end{lem}

The proof of Lemma \ref{lem:1} is given in the Appendix.
It is noted that $\bthh_a^{EB}-\bth_a=(\bthh_a(\bPsi)-\bth_a) + (\bthh_a(\bPsi)-\bthh_a^{EB})$.
From Lemma \ref{lem:1}, $\bthh_a(\bPsi)-\bthh_a^{EB}$ is a function of $\P \y$ and is independent of $\bthh_a(\bPsi)-\bth_a$.
It is noted that 
\begin{align*}
E[(&\bthh_a(\bPsi)-\bth_a)(\bthh_a(\bPsi)-\bth_a)^\top]\\
=&
E[(\btht_a(\bbe, \bPsi)-\bth_a)(\btht_a(\bbe, \bPsi)-\bth_a)^\top + (\bthh_a(\bPsi)-\btht_a(\bbe, \bPsi))(\bthh_a(\bPsi)-\btht_a(\bbe, \bPsi))^\top]
\\
=&\G_{1a}(\bPsi)+\G_{2a}(\bPsi),
\end{align*}
where 
\begin{equation}
\begin{split}
\G_{1a}(\bPsi)=&
(\bPsi^{-1}+\D_a^{-1})^{-1}=\bPsi(\bPsi+\D_a)^{-1}\D_a,\\
\G_{2a}(\bPsi)=&
\D_a(\bPsi+\D_a)^{-1}\X_a\{\X^\top (\I_m\otimes \bPsi+\D)^{-1}\X\}^{-1}\X_a^\top(\bPsi+\D_a)^{-1}\D_a.
\end{split}
\label{eqn:G12}
\end{equation}
Because $\bthh_a(\bPsi)-\bth_a$ is independent of $\P\y$, it is observed that given $\P\y$, the conditional distribution of $\bthh_a(\bPsi)-\bth_a$ is $\Nc_k(\zero, \G_{1a}(\bPsi)+\G_{2a}(\bPsi))$.
This implies the following lemma which will be used for constructing a confidence region.

\begin{lem}
\label{lem:cond}
Under the conditions {\rm (H1)} and {\rm (H2)}, the conditional distribution of $\bthh_a^{EB}-\bth_a$ given $\P  \y$ is
\begin{align}
\bthh_a^{EB}-\bth_a | \P  \y \sim \Nc_k (\bthh_a^{EB}-\bthh_a(\bPsi), \H_a(\bPsi)).
\label{eqn:cd}
\end{align}
for $\H_a(\bPsi)=\G_{1a}(\bPsi)+\G_{2a}(\bPsi)$.
\end{lem}

For evaluating uncertainty of the EBLUP asymptotically, we assume the conditions given below for $m\to\infty$:.

\medskip
(H3)\ $\bPsih$ is $\sqrt{\mathstrut m}$-consistent and second-order unbiased, namely $\bPsih-\bPsi=O(m^{-1/2})$ and $E[\bPsih]=\bPsi+o(m^{-1})$.

\smallskip
(H4)\ $0<k<\infty$, $0<s<\infty$.

\smallskip
(H5)\ There exist positive constants ${\underline d}$ and ${\overline d}$ such that ${\underline d}$ and ${\overline d}$ do not depend on $m$ and satify ${\underline d}\I_k \leq \D_i \leq {\overline d}\I_k$ for $i=1, \ldots, m$.

\smallskip
(H6)\ $\X^\top \X$ is nonsingular and $\X^\top \X/m$ converges to a  positive definite matrix.

\medskip
Under these conditions, we can obtain the important approximations which will be useful for evaluating the mean squared error (MSE) matrix of $\bthh_a^{EB}$ and for constructing corrected confidence region based on $\bthh_a^{EB}$.

\begin{lem}
\label{lem:msemest}
Under conditions {\rm (H1)}-{\rm (H6)}, the following approximations hold:

{\rm (1)}\ $E[\{\bthh_a^{EB}-\bthh_a(\bPsi)\}\{\bthh_a^{EB}-\bthh_a(\bPsi)\}^\top]=\G_{3a}(\bPsi)+O(m^{-3/2})$, where
\begin{align}
\G_{3a}(\bPsi)=\D_a(\bPsi+\D_a)^{-1}E\Big[(\bPsih-\bPsi)(\bPsi+\D_a)^{-1}(\bPsih-\bPsi)\Big](\bPsi+\D_a)^{-1}\D_a,
\label{eqn:G30}
\end{align}

{\rm (2)}\ $E[\G_{1a}(\bPsih)]=\G_{1a}(\bPsi)-\G_{3a}(\bPsi)+O(m^{-3/2})$.
\end{lem}

The proof of Lemma \ref{lem:msemest} is given in the Appendix.
Using Lemma \ref{lem:cond} and Lemma \ref{lem:msemest} (1), we can approximate the MSE matrix of $\bthh_a^{EB}$ as
\begin{align}
{\rm MSEM}(\bthh_a^{EB})=& \G_{1a}(\bPsi)+\G_{2a}(\bPsi) + E[\{\bthh_a^{EB}-\bthh_a(\bPsi)\}\{\bthh_a^{EB}-\bthh_a(\bPsi)\}^\top]
\non\\
=& \G_{1a}(\bPsi)+\G_{2a}(\bPsi) +\G_{3a}(\bPsi) + O(m^{-3/2}).
\label{eqn:MSE}
\end{align}
Using Lemma \ref{lem:msemest} (2), we can obtain the second-order unbiased estimator of ${\rm MSEM}(\bthh_a^{EB})$, which is given by
\begin{equation}
msem(\bthh_a^{EB})=\G_{1a}(\bPsih)+\G_{2a}(\bPsih) +2\G_{3a}(\bPsih),
\label{eqn:msem}
\end{equation}
namely, $E[msem(\bthh_a^{EB})]={\rm MSEM}(\bthh_a^{EB})+O(m^{-3/2})$.
Lemma \ref{lem:msemest} will be also used for deriving corrected confidence region in the next section.

\section{Confidence Region with Corrected Coverage Probability}
\label{sec:cr}

We now construct a confidence region of $\bth_a$ based on $\bthh_a^{EB}$ with second-order accuracy.
When $\bPsi$ is known, it follows from Lemma \ref{lem:cond} that the confidence region based on the  Mahalanobis distance with $100(1-\al)\%$ confidence coefficient is $\{ \bth_a \mid (\bth_a-\bthh_a(\bPsi))^\top \H_a^{-1}(\bPsi) (\bth_a-\bthh_a(\bPsi)) \leq \chi_{k, 1-\al}^2\}$ for $\H_a(\bPsi)=\G_{1a}(\bPsi)+\G_{2a}(\bPsi)$, where $\chi_{k,1-\al}^2$ is the $100\al\%$ upper quantile of the chi-squared distribution with degrees of freedom $k$.
For a matrix $\A(\bPsi)$, $\A^{-1}(\bPsi)$ denotes the inverse matrix of $\A(\bPsi)$. 
Since $\bPsi$ is unknown, we replace $\bPsi$ with estimator $\bPsih$ to get the naive confidence region
\begin{equation}
CR_0=\{ \bth_a \mid (\bth_a-\bthh_a(\bPsih))^\top \H_a^{-1}(\bPsih) (\bth_a-\bthh_a(\bPsih)) \leq \chi_{k, 1-\al}^2\}.
\label{eqn:conf0}
\end{equation}
Under appropriate conditions, it can be shown that the coverage probability tends to the nominal confidence coefficient $1-\al$, namely $\lim_{m\to\infty}P(\bth_a \in CR_0)=1-\al$.
However, this confidence region has the second-order bias, because $P(\bth_a \in CR_0)=1-\al + O(m^{-1})$.
Thus, we want to derive a corrected confidence region $CR$ such that $P(\bth_a\in CR)=1-\al+O(m^{-3/2})$.

\medskip
Define $B_1$, $B_2$ and $B_3$ by
\begin{equation}
\begin{split}
B_1=& B_1(\bPsi)=
-{1 \over 2}\tr\Big( E[ \K_{a}(\bPsih) \H_a^{-1}(\bPsi) \K_{a}(\bPsih))]\Big),
\\
B_2=&B_2(\bPsi)= -{1 \over 8}\Big\{E[\tr^2( \K_{a}(\bPsih))]+2\tr\Big( E[( \K_{a}(\bPsih))^2]\Big)\Big\},
\\
B_3=&\tr(\H_a^{-1}(\bPsi)\G_{3a}(\bPsi)),
\end{split}
\label{eqn:B12}
\end{equation}
where $\K_{a}(\bPsih)=\H_a^{-1/2}(\bPsi)(\G_{1a}(\bPsih)-\G_{1a}(\bPsi))\H_a^{-1/2}(\bPsi)$ and $\tr^2(\A)=(\tr \A)^2$ for matrix $\A$.
It can be seen that $B_1=O(m^{-1})$, $B_2=O(m^{-1})$ and $B_3=O(m^{-1})$.
Then, we provide the main theorem which will be proved in the Appendix.

\begin{thm}\label{thm:cr}
Under the conditions {\rm (H1)}-{\rm (H6)}, it holds that
\begin{align}
P\{(\bthh_a^{EB}&-\bth_a)^\top\H_a^{-1}(\bPsih)(\bthh_a^{EB}-\bth_a)\leq x\}
\non\\
&=F_{k}(x)+2(B_1-B_3-B_2)f_{k+2}(x)+2B_2f_{k+4}(x)+o(m^{-1}),
\label{eqn:main}
\end{align}
where $F_{k}(x)$ and $f_k(x)$ are the cumulative distribution and probability density functions of the chi-squared distribution with the degree of freedom $k$, respectively.
\end{thm}

We can consider the Bartlett-type correction using the asymptotic expansion (\ref{eqn:main}).
For $h=O(m^{-1})$, it is observed that
\begin{align*}
P\{(\bthh_a^{EB}&-\bth_a)^\top\H_a^{-1}(\bPsih)(\bthh_a^{EB}-\bth_a)\leq x(1+h)\}
\\
=&
F_{k}(x)+hxf_k(x)+2(B_1-B_3-B_2)f_{k+2}(x)+2B_2f_{k+4}(x)+o(m^{-1}).
\end{align*}
Note that $hxf_k(x)+2(B_1-B_3-B_2)f_{k+2}(x)+2B_2f_{k+4}(x)$ is of order $O(m^{-1})$.
Thus, the second-order term vanishes if 
\begin{equation}
hx f_k(x)=-2(B_1-B_3-B_2)f_{k+2}(x)-2B_2f_{k+4}(x)=0.
\label{eqn:bart}
\end{equation}
Since $\Ga(x+1)=x\Ga(x)$ for the gamma function $\Ga(x)$, the solution of the equation (\ref{eqn:bart}) on $h$ is 
\begin{equation}
h^*(\bPsi)=-2\{(B_1-B_3-B_2)/k+B_2x/k(k+2)\}.
\label{eqn:h}
\end{equation}
For $h^*(\bPsi)$ given in $(\ref{eqn:h})$, it holds that for any $x>0$,
\begin{align*}
P\{(1+h^*(\bPsih))^{-1}(\bthh_a^{EB}-\bth_a)^\top\H_a^{-1}(\bPsih)(\bthh_a^{EB}-\bth_a)\leq x\}=F_{k}(x)+o(m^{-1}).
\end{align*}
Hence, the corrected confidence region is given by
\begin{equation}
CR=\{ \bth_a \mid (\bth_a-\bthh_a^{EB})^\top \H_a^{-1}(\bPsih) (\bth_a-\bthh_a^{EB}) \leq \{1+h^*(\bPsih)\}\chi_{k, 1-\al}^2\}.
\label{eqn:conf}
\end{equation}

\begin{cor}
Under conditions {\rm (H1)}-{\rm (H6)}, it holds that
$$
P(\bth_a \in CR) = 1-\al + o(m^{-1}).
$$
\end{cor}

\section{Derivation of a Second-order Unbiased and Positive-definite Estimator of $\bPsi$}
\label{sec:psi}

We here provide a new method for deriving a second-order unbiased and positive-definite estimator of $\bPsi$.
As well known in the univariate case, the Prasad-Rao estimator of the ^^ between' component of variance takes a negative value with a positive probability, and the nonnegative estimator which truncates it at zero is used.
The maximum likelihood (ML) and restricted maximum likelihood (REML) estimators take values of zero with positive probabilities.
To fix this drawback, Li and Lahiri (2010) suggested the adusted maximum likelihood method for giving a positive and consistent estimator.
As pointed out by Yoshimori and Lahiri (2014), this problem causes instability of the corrected confidence interval.
In the multivariate case, since $\G_{2a}(\bPsi)=O(m^{-1})$, it is seen that $\H_a^{-1}(\bPsi)=\G_{1a}^{-1}(\bPsi)+O(m^{-1})=\bPsi^{-1}+\D_a^{-1}+O(m^{-1})$.
This means that the correction function $h^*(\bPsi)$ takes a large value when some eigenvalues of estimator $\bPsih$ are zero.

\medskip
To derive a positive-definite and consistet estimator of $\bPsi$, let $\U$ be a $k\times k$ orthogonal matrix $\U$ such that $\bPsih=\U\text{\boldmath $L$} \U^\top$ for a diagonal matrix $\text{\boldmath $L$}=\diag(\ell_1, \ldots, \ell_k)$.
Then, we consider adjusted estimators of the form
\begin{align}
\bPsih_{(A)} = {1 \over 2}(\bPsih-a\I_k+\U\text{\boldmath $L$}_{(A)} \U^\top),
\label{eqn:Psis}
\end{align}
where
$$
\text{\boldmath $L$}_{(A)}=\diag(\sqrt{\mathstrut (\ell_1-\ah)^2+\bh_1}, \ldots, \sqrt{\mathstrut (\ell_k-\ah)^2+\bh_k}),
$$
for some statistics $\ah$ and $\bh_1, \ldots, \bh_k$.

\begin{prp}\label{prp:pd}
Assume that $\ah=O_p(m^{-1})$, $E[\ah^2]=o(m^{-1})$ and that $\bh_i$'s are positive almost surely and $\bh_i=O_p(m^{-1})$ for $i=1,\ldots,k$.
Let $\bPsih$ be a consistent estimator of $\bPsi$ as $m\to\infty$.

{\rm (1)}\ $\bPsih_{(A)}$ given in $(\ref{eqn:Psis})$ is positive-definite almost surely, and $\bPsih_{(A)}=\bPsih+O_p(m^{-1})$.

{\rm (2)}\ If $\bPsih$ is second-order unbiased, namely, $E[\bPsih]=\bPsi+o(m^{-1})$, and if  $\bh_i=4\ah(\ell_i-\ah)0$ is almost surely positive, then $\bPsih_{(A)}$ is positive definite almost surely and second-order unbiased. 
\end{prp}

{\bf Proof}.\ \ 
It is clear that $\bPsih_{(A)}$ is positive definite almost surely.
Note that there exists positive $\la_i$ such that $\ell_i$ converges to $\la_i$, because $\bPsih$ is consistent.
Since $\ah=O_p(m^{-1})$ and $E[\ah^2]=o(m^{-1})$, it is seen that $P(\la_i-\ah <0)=o(m^{-1})$.
Then, the eigenvalues of $\bPsih_{(A)}$ are approximated as
\begin{align}
\ell_i-\ah+\sqrt{\mathstrut (\ell_i-\ah)^2+\bh_i}=&\ell_i-\ah+\sqrt{\mathstrut (\ell_i-\la_i+\la_i-\ah)^2+\bh_i}
\non\\
=&
\ell_i-\ah+|\la_i-\ah|\sqrt{\mathstrut1+{2(\la_i-\ah)(\ell_i-\la_i)+(\ell_i-\la_i)^2+\bh_i \over  (\la_i-\ah)^2}}
\non\\
=&
2\ell_i-2\ah+{\bh_i \over 2(\la_i-\ah)}+o_p(m^{-1}).
\label{eqn:prp1}
\end{align}
This implies that $\bPsih_{(A)}=\bPsih+O_p(m^{-1})$, which shows part (1).
For part (2), let $\bh_i=4\ah(\la_i-\ah)$.
Then we can see that the second term is equal to the third term in RHS of (\ref{eqn:prp1}), and the second-order bias vanishes.
Thus, the part (2) is shown by replacing $\la_i$ with $\ell_i$.
\hfill$\Box$

\bigskip
Before constructing the estimator $\bPsih_{(A)}$ with specific $\ah$ and $\bh_i$'s, we obtain estimator $\bPsih$ which satisfies conditions (H1), (H2) and (H3).
When $\bPsi$ is a fully unknown covariance matrix, it is hard to derive the ML and REML estimates numerically.
Instead, we begin by deriving a Prasad-Rao type estimator based on the moment method.
Because $E[(\y_i-\X_i\bbe)(\y_i-\X_i\bbe)^\top]=\bPsi+\D_i$ for $i=1, \ldots, m$, we have $\sum_{i=1}^m E[(\y_i-\X_i\bbe)(\y_i-\X_i\bbe)^\top]=m \bPsi+\sum_{i=1}^m\D_i$.
Substituting the ordinary least squares estimator $\bbeh^{OLS}=(\X^\top\X)^{-1}\X^\top\y$ into $\bbe$, we get the Prasd-Rao type consistent estimator
\begin{equation}
\bPsih_0^{PR} = {1\over m} \sum_{i=1}^m \big\{ (\y_i-\X_i\bbeh^{OLS})(\y_i-\X_i\bbeh^{OLS})^\top-\D_i\big\}.
\label{eqn:Psi0}
\end{equation}
It is noted that this estimator has a second-order bias.
In fact, the bias, given by ${\rm Bias}_{\bPsih_0^{PR}}(\bPsi)=E[\bPsih_0^{PR}]-\bPsi$, is 
\begin{align}
{\rm Bias}_{\bPsih_0^{PR}}(\bPsi)=&
 {1\over m}\sum_{i=1}^m \X_i(\X^\top\X)^{-1}\Big\{\sum_{j=1}^m\X_j^\top (\bPsi+\D_j)\X_j\Big\}(\X^\top\X)^{-1} \X_i^\top
\non\\
&
-{1\over m}\sum_{i=1}^m (\bPsi+\D_i)\X_i(\X^\top\X)^{-1}\X_i^\top
-{1\over m}\sum_{i=1}^m \X_i(\X^\top\X)^{-1}\X_i^\top(\bPsi+\D_i).
\label{eqn:Bias}
\end{align}
Substituting $\bPsih_0^{PR}$ into ${\rm Bias}_{\bPsih_0^{PR}}(\bPsi)$ provides the bias-corrected estimator
\begin{equation}
\bPsih^{PR} = \bPsih_0^{PR} - {\rm Bias}_{\bPsih_0^{PR}}(\bPsih_0^{PR}).
\label{eqn:Psi1}
\end{equation}

The estimator $\bPsih^{PR}$ satisfies conditions (H1), (H2) and (H3).
However, it still has a drawback of taking a negative value with a positive probability.
For applying the method suggested in Proposition \ref{prp:pd}, let 
$$
\ah=\tr(\bPsih^{PR})/mk\quad {\rm and}\quad \bh_i=\max\{4\ah(\ell_i^{PR}-\ah), 1/m\}, \quad
{\rm for}\quad i=1,\ldots,k,
$$
where $\ell_i^{PR}$'s are eigenvalues of $\bPsih^{PR}$.
Note that $P\{\ah(\ell_i^{PR}-\ah)<1/(4m)\}=o(m^{-1})$.
Then, we suggest the adjusted estimator
\begin{align}
\bPsih_{(A)}^{PR} = {1 \over 2}(\bPsih^{PR}-\ah\I_k+\U^{PR}\text{\boldmath $L$}^{PR}_{(A)} (\U^{PR})^\top),
\label{eqn:Psia}
\end{align}
where column vectors of $\U^{PR}$ are the eigenvectors of $\bPsih^{PR}$ and
$$
\text{\boldmath $L$}_{(A)}^{PR}=\diag(\sqrt{\mathstrut (\ell_1^{PR}-\ah)^2+\bh_1}, \ldots, \sqrt{\mathstrut (\ell_k^{PR}-\ah)^2+\bh_k}).
$$
It follows from Proposition \ref{prp:pd} that $\bPsih_{(A)}^{PR}$ is positive-definite and second-order unbiased.

\medskip
Before calculating some moments given in $B_1$, $B_2$ and $B_3$, we need a closed-form expression of $\G_{3a}(\bPsi)$ given in (\ref{eqn:G30}), which is stated in the following lemma.
 
\begin{lem}
\label{lem:MSE}
By using the Prasad-Rao type estimetor given in $(\ref{eqn:Psi1})$ or $(\ref{eqn:Psia})$, we can write $\G_{3a}(\bPsi)$ in $(\ref{eqn:G30})$, as
\begin{align}
\begin{split}
\G_{3a}(\bPsi)=&
{1\over m^2}\D_a(\bPsi+\D_a)^{-1}
\Big[ \sum_{i=1}^m(\bPsi+\D_i)(\bPsi+\D_a)^{-1}(\bPsi+\D_i) 
\\
&\qquad+ \sum_{i=1}^m\{\tr[(\bPsi+\D_i) (\bPsi+\D_a)^{-1}]\}(\bPsi+\D_i)\Big](\bPsi+\D_a)^{-1}\D_a.
\end{split}
\label{eqn:G3}
\end{align}
\end{lem}

\medskip
Finally, we calculate some moments given in $B_1$, $B_2$ and $B_3$ for the estimator $\bPsih_{(A)}^{PR}$.
This calculation is used for providing the correction function $h^*(\bPsi)$.

\begin{lem}\label{lem:moment}
Assume conditions {\rm (H4)}-{\rm (H6)}.
For $\bPsih_{(A)}^{PR}$ as in $(\ref{eqn:Psia})$, the values of $B_1$ and $B_2$ in $(\ref{eqn:B12})$ are given by
\begin{align*}
B_1
=&
-{1 \over 2m^2} \sum_{i=1}^m \Big\{\tr\Big((\bPsi+\D_a)^{-1}\D_a\H_a^{-2}(\bPsi)\D_a(\bPsi+\D_a)^{-1} (\bPsi+\D_i)
\\
&\qquad\qquad\times(\bPsi+\D_a)^{-1}\D_a\H_a^{-1}(\bPsi)\D_a(\bPsi+\D_a)^{-1}(\bPsi+\D_i)\Big)
\\
&\qquad+\tr\Big((\bPsi+\D_a)^{-1}\D_a\H_a^{-2}(\bPsi)\D_a(\bPsi+\D_a)^{-1} (\bPsi+\D_i)\Big)
\\
&\qquad\qquad\times\tr\big((\bPsi+\D_a)^{-1}\D_a\H_a^{-1}(\bPsi)\D_a(\bPsi+\D_a)^{-1}(\bPsi+\D_i)\Big)\Big\}
+o(m^{-1}),
\\
B_2=&
-{ 1\over 4m^2}\sum_{i=1}^{m} \Big\{ \tr\Big(((\bPsi+\D_a)^{-1}\D_a\H_a^{-1}(\bPsi)\D_a(\bPsi+\D_a)^{-1} (\bPsi+\D_i))^2\Big)
\\
&\qquad+\tr\Big(((\bPsi+\D_a)^{-1}\D_a\H_a^{-1}(\bPsi)\D_a(\bPsi+\D_a)^{-1}(\bPsi+\D_i))^2\Big)
\\
&\qquad+\tr^2\Big(((\bPsi+\D_a)^{-1}\D_a\H_a^{-1}(\bPsi)\D_a(\bPsi+\D_a)^{-1} (\bPsi+\D_i))\Big)\Big\}
+o(m^{-1}),
\end{align*}
and the value of $B_3$ in $(\ref{eqn:B12})$ is $B_3=\tr(\H_a^{-1}(\bPsi)\G_{3a}(\bPsi))$ for $\G_{3a}(\bPsi)$ given in $(\ref{eqn:G3})$.
\end{lem}

By substituting these values into (\ref{eqn:h}), we can construct the confidence region in the closed-form.
Moreover, by substituting (\ref{eqn:G3}) into (\ref{eqn:msem}), we can obtain an estimator of closed-form approximation of the MSE matrix of $\bthh_a^{EB}$ as a by-product.

\section{Confidence Region for the Difference of Two Small Area Means}
\label{sec:crd}

In this section, we extend the results in Section \ref{sec:cr} to the construction of a confidence region for $\bth_a-\bth_b$ for $a\neq b$.
This enables us to conduct a statistical test under the null hypothsis $H_0:\bth_a=\bth_b$.
Since the corrected confidence region of $\bth_a-\bth_b$ can be constructed by the same arguments as in Section \ref{sec:cr}, we here provide the sketch of the result.

\medskip
Let $\G_{ab}(\bPsi)=\H_a(\bPsi)+\H_b(\bPsi)-\G_{2ab}(\bPsi)$, where $\G_{2ab}(\bPsi)=E[(\bthh_a^{EB}-\bth_a)(\bthh_b^{EB}-\bth_b)^\top]+E[(\bthh_b^{EB}-\bth_b)(\bthh_a^{EB}-\bth_a)^\top]$.
Then, it can be evaluated as
\begin{align*}
\G_{2ab}(\bPsi)=&
\D_a(\bPsi+\D_a)^{-1}\X_a\{\X^\top (\I_m\otimes \bPsi+\D)^{-1}\X\}^{-1}\X_b^\top(\bPsi+\D_b)^{-1}\D_b
\\
&+\D_b(\bPsi+\D_b)^{-1}\X_b\{\X^\top (\I_m\otimes \bPsi+\D)^{-1}\X\}^{-1}\X_a^\top(\bPsi+\D_a)^{-1}\D_a.
\end{align*}
The asymptotic expansion of the cumulative distribution function is
\begin{align*}
P\{&(\bthh_a^{EB}-\bth_a-\bthh_b^{EB}+\bth_b)^\top\G_{ab}^{-1}(\bPsih)(\bthh_a^{EB}-\bth_a-\bthh_b^{EB}+\bth_b)\leq x\}\\
=&F_{k}(x)+2({\tilde B_1} - {\tilde B_3} - {\tilde B_2})f_{k+2}(x)+2\tilde B_2f_{k+4}(x)+o(m^{-1}),
\end{align*}
where $\tilde B_1$, $\tilde B_2$ and ${\tilde B_3}$ are 
\begin{equation}
\begin{split}
\tilde B_1=&
-{1 \over 2}\tr\Big( E[\G_{ab}^{-1/2}(\bPsi)(\G_{1a}(\bPsih)-\G_{1a}(\bPsi)+\G_{1b}(\bPsih)-\G_{1b}(\bPsi))\G_{ab}^{-2}(\bPsi)
\\
&\qquad\qquad\times(\G_{1a}(\bPsih)-\G_{1a}(\bPsi)+\G_{1b}(\bPsih)-\G_{1b}(\bPsi))\G_{ab}^{-1/2}(\bPsi)]\Big),
\\
\tilde B_2=& -{1 \over 8}\Big\{E[\tr^2( (\G_{1a}(\bPsih)-\G_{1a}(\bPsi)+\G_{1b}(\bPsih)-\G_{1b}(\bPsi))\G_{ab}^{-1}(\bPsi))]
\\
&\qquad+2\tr\Big( E[((\G_{1a}(\bPsih)-\G_{1a}(\bPsi)+\G_{1b}(\bPsih)-\G_{1b}(\bPsi))\G_{ab}^{-1}(\bPsi))^2]\Big)\Big\},
\\
\tilde B_3=&\tr(\G_{ab}^{-1}(\bPsi)(\G_{3a}(\bPsi)+\G_{3b}(\bPsi)-\G_{2ab}(\bPsi))),
\end{split}
\label{eqn:tB12}
\end{equation}
Setting $\tilde h^*=-2\{({\tilde B_1}-{\tilde B_3}-{\tilde B_2})/k+\tilde B_2x/k(k+2)\}$, we have
\begin{align*}
P((1+\tilde h^*)^{-1}(\bthh_a^{EB}-\bth_a-\bthh_b^{EB}+\bth_b)^\top\G_{ab}^{-1}(\bPsih)(\bthh_a^{EB}-\bth_a-\bthh_b^{EB}+\bth_b)\leq x)=F_{k}(x)+o(m^{-1}),
\end{align*}
namely, $P(\bth_a-\bth_b \in CR_{ab})=1-\al+o(m^{-1})$ for the corrected confidence region
$$
CR_{ab}=\{\bth_a-\bth_b | 
(\bthh_a^{EB}-\bth_a-\bthh_b^{EB}+\bth_b)^\top\G_{ab}^{-1}(\bPsih)(\bthh_a^{EB}-\bth_a-\bthh_b^{EB}+\bth_b)\leq (1+\tilde h^*)\chi_{k, 1-\al}^2\}.
$$

\medskip
When the adjusted Prasd-Rao type estimator $\bPsih_{(A)}^{PR}$ given in (\ref{eqn:Psia}) is used for estimating $\bPsi$, the functions $\tilde B_1$ and $\tilde B_2$ are calculated as
Then, we have
\begin{align*}
{\tilde B_1}=&
-{1\over 2}\tr(\V_{1aa}+\V_{1bb}+\V_{1ab}+\V_{1ba}) +o(m^{-1}),
\\
{\tilde B_2}=&
- {1 \over 4m^2}\sum_{i=1}^{m}\tr\Big(\Big\{(\bPsi+\D_a)^{-1}\D_a\G_{ab}^{-1}\D_a(\bPsi+\D_a)^{-1} (\bPsi+\D_i)\\
&\qquad
+(\bPsi+\D_b)^{-1}\D_b\G_{ab}^{-1}\D_b(\bPsi+\D_b)^{-1} (\bPsi+\D_i)\Big\}^2\Big)
\\
& -{1\over 4} \tr(\V_{2aa}+\V_{2bb}+\V_{2ab}+\V_{2ba})+o(m^{-1}),
\end{align*}
where for $(c,d)=(a,a), (a,b), (b,a)$ and $(b,b)$,
\begin{align*}
\V_{1cd}
=&
\G_{cd}^{-1/2}(\bPsi)\D_c(\bPsi+\D_c)^{-1}\\
&\times \Big[{1 \over m^2} \sum_{i=1}^m \Big\{ (\bPsi+\D_i)(\bPsi+\D_c)^{-1}\D_c\G_{cd}^{-2}(\bPsi)\D_d(\bPsi+\D_d)^{-1} (\bPsi+\D_i)
\\
&\qquad
+ \tr((\bPsi+\D_c)^{-1}\D_c\G_{cd}^{-2}(\bPsi)\D_d(\bPsi+\D_d)^{-1} (\bPsi+\D_i))(\bPsi+\D_i)\Big\}\Big]\\
&\times
\D_d(\bPsi+\D_d)^{-1}\G_{cd}^{-1/2}(\bPsi),
\\
\V_{2cd}
=&
\G_{cd}^{-1/2}(\bPsi)
\D_c(\bPsi+\D_c)^{-1} 
\\
&\times\Big[{1 \over m^2} \sum_{i=1}^m \Big\{ (\bPsi+\D_i)(\bPsi+\D_c)^{-1}\D_c\G_{cd}^{-1}(\bPsi)\D_d(\bPsi+\D_d)^{-1} (\bPsi+\D_i)
\\
&\qquad
+ \tr((\bPsi+\D_c)^{-1}\D_c\G_{cd}^{-1}(\bPsi)\D_d(\bPsi+\D_d)^{-1} (\bPsi+\D_i))(\bPsi+\D_i)\Big\}\Big]
\\
&\times (\bPsi+\D_d)^{-1}\D_d\G_{cd}^{-1/2}(\bPsi).
\end{align*}
Also, $\tilde B_3$ is obtained by using the expression in (\ref{eqn:G3}).

\section{Finite Sample Performances}
\label{sec:sim}

We now investigate finite sample performances of the proposed confidence regions by simulation in  the multivariate Fay-Herriot model (\ref{eqn:MFH}) for $k=2, 3$ and $m=30$.
The design matrix, $\X_{i}$ is a $k \times 2k$ matrix, such that 
$$
\X_{i}=\begin{pmatrix} 1 & x_{i1} & 0 & 0\\ 0&0&1 & x_{i2} \end{pmatrix}, \X_{i}=\begin{pmatrix} 1 & x_{i1} & 0 & 0 & 0 & 0\\ 0&0&1 & x_{i2} & 0 & 0\\ 0&0&0&0& 1& x_{i3} \end{pmatrix}
$$
for $k=2,3$ respectevely, where $x_{ij}$'s are generated from the uniform distribution on $(-1,1)$, which are fixed through the simulation runs.
As a setup of the covariance matrix $\bPsi$ of the random effects, we consider 
$$
\bPsi= \left\{\begin{array}{ll}
\rho\bpsi_2\bpsi_2^\top+(1-\rho){\rm diag}(\bpsi_2\bpsi_2^\top) & {\rm for}\ k=2,\\
\rho\bpsi_3\bpsi_3^\top+(1-\rho){\rm diag}(\bpsi_3\bpsi_3^\top) & {\rm for}\ k=3,
\end{array}\right.
$$
where $\bpsi_2=(\sqrt{1.6}, \sqrt{0.8})^\top$, $\bpsi_3=(\sqrt{1.6}, \sqrt{1.2}, \sqrt{0.8})^\top$, and ${\rm diag}(\A)$ denotes the diagonal matrix consisting of diagonal elements of matrix $\A$.
Here, $\rho$ is the correlation coefficient, and  we handle the three cases $\rho=0.2, 0.4, 0.6$.
The cases of negative correlations are omitted, because we observe the same results with those of positive ones.

\medskip
Concerning the dispersion matrices $\D_i$ of sampling errors $\bep_i$, we treat the two $\D_i$-patterns: (a) $0.7{\bf I}_k$, $0.6{\bf I}_k$, $0.5{\bf I}_k$, $0.4{\bf I}_k$, $0.3{\bf I}_k$ and (b) $2.0{\bf I}_k$, $0.6{\bf I}_k$, $0.5{\bf I}_k$, $0.4{\bf I}_k$, $0.2{\bf I}_k$.
These cases are treated by Datta, $\et$ (2005) in the univariate Fay-Herriot model. 
There are five groups $G_1, \ldots, G_5$ corresponding to these $\D_i$-patterns, 
 and there are six small areas in each group for $m=30$, respectively, where the sampling covariance matrices $\D_i$ are the same for areas within the same group.
Concerning the underlying distributions for $\v_i$ and $\bep_i$, we consider two kinds of distributions, that is, multivariate normal distributions and multivariate normalized chi-squared distributions with degrees of freedom $2$, which are denoted by M1 and M2, respectively.
The chis-quared distribution is used for investigating robustness of the proposed method against the misspecification of distributions of $\v_i$ and $\bep_i$.
The values of coverage probabilities (CP) of the corrected confidence region and the naive confidence region and the values of the Bartlett-type correction term $h^*$ are obtained based on $10,000$ simulation run, where the nominal confidence coefficient is $95\%$.

\medskip
The values of CP and the correction term in the case of $k=2$ are reported in Tables \ref{tab:cia} and \ref{tab:cib} for $\D_i$-patterns (a) and (b), respectively.
From the tables for normal distributions, the corrected method has CP values larger than the nominal confidence coefficient.
In contrast, CP values of the naive confidence region are much smaller than the nominal confidence coefficient.
For example, CP value for $G_1$ in Table \ref{tab:cib} is about $89\%$.
These show that the naive method is not appropriate for a confidence region and the correction by $h^*$ works well.
For chi-square distributions, CP values of the corrected method satisfies the nominal confidence coefficient in most cases except few cases where CP values are slightly smaller than, but close to $95\%$, while the performance of the naive method is worse than that in the normal distributions.
Thus, the corrected method remains good and robust for the chi-square distributions.
Concerning the Bartlett-type correction, it increases as sampling variances or correlation coefficients $\rho$ increase.

\medskip
Table \ref{tab:cik3} reports the results for $k=3$ and $\D_i$-pattern (a).
Comparing Tables \ref{tab:cia} and \ref{tab:cik3}, we can observe that CP values of the naive confidence region are worse in $k=3$ than those in $k=2$.
The corrected confidence region satisfies the nominal confidence coefficient for $k=3$ in most cases except the case of $\rho=0.2$ in chi-square distributions.
Hence, the corrected method works well and is robust still for $k=3$.

\medskip
We next investigate the finite sample performance of the corrected confidence region for the difference of two small area means, $\bth_a-\bth_b$ for $k=2$ and $\D_i$-pattern (a), where the corrected method is provided in Section \ref{sec:crd}.
In each area group, we consider the difference between the first two small areas means.
Table \ref{tab:cida} reports values of the coverage provabilities (CP) and the Bartlett-type correction term $h^*$ for $\bth_a-\bth_b$.
From Table \ref{tab:cida}, it is revealed that the performances are similar to the results in Table \ref{tab:cia}, while values of the Bartlett-type correction term $h^*$ are larger for $\bth_a-\bth_b$.

\begin{table}[htbp]
\begin{center}
\resizebox{12cm}{!} {
\begin{tabular}{cccccccc}
 \hline
&&\multicolumn{3}{c}{Normal}&\multicolumn{3}{c}{chi-square} \\
\cline{3-5}\cline{6-8}
$\rho$& & 0.2 & 0.4 & 0.6 & 0.2 & 0.4 & 0.6   \\
\hline
$G_1$& CP & 0.955 & 0.968 & 0.974 & 0.939 & 0.945 & 0.956 \\ 
&            & (0.917) & (0.923) & (0.917) & (0.898) & (0.901) & (0.907) \\ 
& $h^*$ & 0.429 & 0.492  & 0.760 & 0.636 & 0.697 & 0.841 \\ 
$G_2$& CP & 0.962 & 0.960 & 0.977 & 0.941 & 0.942 & 0.954  \\  
&               & (0.923) & (0.913) & (0.922) & (0.902) & (0.899) & (0.912) \\
& $h^*$ & 0.534 & 0.571 & 0.865 & 0.598 & 0.640 &  0.758 \\ 
$G_3$& CP & 0.958 & 0.962 & 0.978 & 0.939 & 0.947 & 0.953 \\ 
&               & (0.921) & (0.921) & (0.922) & (0.901) & (0.906) & (0.912) \\  
& $h^*$ & 0.470 & 0.530 & 0.843 & 0.669 & 0.731 & 0.849 \\ 
$G_4$& CP & 0.959 & 0.965 & 0.973 & 0.939 & 0.944 & 0.953\\  
&               & (0.928) & (0.928) & (0.925) & (0.905) & (0.908) & (0.911) \\ 
& $h^*$ & 0.388 & 0.441 & 0.688 & 0.552 & 0.610 &  0.742\\ 
$G_5$& CP & 0.954 & 0.962 & 0.976 & 0.951 & 0.947 & 0.955 \\  
&               & (0.923) & (0.927) & (0.930) & (0.914) & (0.914) & (0.924) \\ 
& $h^*$ & 0.441 & 0.470 & 0.734 & 0.480 & 0.519 & 0.622 \\ 
\hline
\end{tabular}
}
\caption{{Coverage probabilities (CP) for nominal $95\%$ confidence regions for $k=2$ and $\D_i$-pattern (a).
(the corrected method in the first line and the naive method in parentheses) }}
\label{tab:cia}
\end{center}
\end{table}

\begin{table}[htbp]
\begin{center}
\resizebox{12cm}{!} {
\begin{tabular}{cccccccc}
 \hline
&&\multicolumn{3}{c}{Normal}&\multicolumn{3}{c}{chi-square} \\
\cline{3-5}\cline{6-8}
$\rho$& & 0.2 & 0.4 & 0.6 & 0.2 & 0.4 & 0.6   \\
\hline
$G_1$& CP & 0.974 & 0.980 & 0.990 & 0.952 & 0.957 & 0.963 \\ 
&            & (0.897) & (0.895) & (0.899) & (0.876) & (0.876) & (0.891) \\ 
& $h^*$ & 1.288 & 1.645  & 2.571 & 1.962 & 2.145 & 2.270 \\ 
$G_2$& CP & 0.969 & 0.980 & 0.987 & 0.953 & 0.960 & 0.967  \\  
&               & (0.905) & (0.908) & (0.909) & (0.885) & (0.894) & (0.897) \\
& $h^*$ & 1.493 & 1.870 & 2.979 & 2.207 & 2.380 &  2.466 \\ 
$G_3$& CP & 0.967 & 0.976 & 0.984 & 0.954 & 0.961 & 0.969 \\ 
&               & (0.912) & (0.908) & (0.908) & (0.894) & (0.898) & (0.901) \\  
& $h^*$ & 1.288 & 1.730 & 2.876 & 1.933 & 2.173 & 2.332 \\ 
$G_4$& CP & 0.967 & 0.974 & 0.982 & 0.953 & 0.958 & 0.965\\  
&               & (0.916) & (0.918) & (0.914) & (0.898) & (0.898) & (0.907) \\ 
& $h^*$ & 1.107 & 1.433 & 2.319 & 1.695 & 1.854 & 1.950 \\ 
$G_5$& CP & 0.966 & 0.973 & 0.980 & 0.954 & 0.957 & 0.965 \\  
&               & (0.926) & (0.921) & (0.922) & (0.904) & (0.909) & (0.917) \\ 
& $h^*$ & 1.169 & 1.494 & 2.471 & 1.696 & 1.857 & 1.928 \\ 
\hline
\end{tabular}
}
\caption{{Coverage probabilities (CP) for nominal $95\%$ confidence regions for $k=2$ and $\D_i$-pattern (b).
(the corrected method in the first line and the naive method in parentheses)}}
\label{tab:cib}
\end{center}
\end{table}

\begin{table}[htbp]
\begin{center}
\resizebox{12cm}{!} {
\begin{tabular}{cccccccc}
 \hline
&&\multicolumn{3}{c}{Normal}&\multicolumn{3}{c}{chi-square} \\
\cline{3-5}\cline{6-8}
$\rho$& & 0.2 & 0.4 & 0.6 & 0.2 & 0.4 & 0.6   \\
\hline
$G_1$& CP & 0.964 & 0.977 & 0.987 & 0.941 & 0.953 & 0.964 \\ 
&            & (0.897) & (0.903) & (0.917) & (0.877) & (0.884) & (0.887) \\ 
& $h^*$ & 0.527 & 0.675  & 0.809 & 0.816 & 0.890 & 1.284 \\ 
$G_2$& CP & 0.964 & 0.976 & 0.989 & 0.950 & 0.951 & 0.965  \\  
&               & (0.897) & (0.897) & (0.920) & (0.883) & (0.878) & (0.886) \\
& $h^*$ & 0.570 & 0.734 & 0.882 & 0.884 & 0.975 &  1.440 \\ 
$G_3$& CP & 0.966 & 0.975 & 0.986 & 0.943 & 0.955 & 0.967 \\ 
&               & (0.903) & (0.903) & (0.917) & (0.879) & (0.884) & (0.889) \\  
& $h^*$ & 0.579 & 0.744 & 0.876 & 0.891 & 0.981 & 1.454 \\ 
$G_4$& CP & 0.965 & 0.973 & 0.985 & 0.940 & 0.952 & 0.965\\  
&               & (0.908) & (0.910) & (0.923) & (0.881) & (0.889) & (0.893) \\ 
& $h^*$ & 0.488 & 0.630 & 0.755 & 0.772 & 0.843 &  1.252\\ 
$G_5$& CP & 0.964 & 0.972 & 0.983 & 0.944 & 0.953 & 0.968 \\  
&               & (0.916) & (0.912) & (0.922) & (0.893) & (0.898) & (0.904) \\ 
& $h^*$ & 0.474 & 0.610 & 0.727 & 0.761 & 0.826 & 1.247 \\ 
\hline
\end{tabular}
}
\caption{{Coverage probabilities (CP) for nominal $95\%$ confidence regions for $k=3$ and $\D_i$-pattern (a).
(the corrected method in the first line and the naive method in parentheses)}}
\label{tab:cik3}
\end{center}
\end{table}

\begin{table}[htbp]
\begin{center}
\resizebox{12cm}{!} {
\begin{tabular}{cccccccc}
 \hline
&&\multicolumn{3}{c}{Normal}&\multicolumn{3}{c}{chi-square} \\
\cline{3-5}\cline{6-8}
$\rho$& & 0.2 & 0.4 & 0.6 & 0.2 & 0.4 & 0.6   \\
\hline
$G_1$& CP & 0.975 & 0.983 & 0.990 & 0.960 & 0.966 & 0.973 \\ 
&            & (0.912) & (0.919) & (0.928) & (0.895) & (0.911) & (0.911) \\ 
& $h^*$ & 0.940 & 0.957  & 1.364 & 1.114 & 1.373 & 1.969 \\ 
$G_2$& CP & 0.977 & 0.986 & 0.986 & 0.960 & 0.967 & 0.975  \\  
&               & (0.922) & (0.933) & (0.929) & (0.907) & (0.912) & (0.927) \\
& $h^*$ & 0.742 & 0.770 & 1.087 & 0.914 & 1.106 & 1.543  \\ 
$G_3$& CP & 0.971 & 0.971 & 0.988 & 0.959 & 0.961 & 0.971 \\ 
&               & (0.928) & (0.917) & (0.932) & (0.912) & (0.908) & (0.911) \\  
& $h^*$ & 0.715 & 0.727 & 1.028 & 0.870 & 1.067 & 1.491 \\ 
$G_4$& CP & 0.974 & 0.976 & 0.983 & 0.950 & 0.965 & 0.978\\  
&               & (0.936) & (0.931) & (0.924) & (0.904) & (0.912) & (0.913) \\ 
& $h^*$ & 0.675 & 0.712 & 1.035 & 0.834 & 1.049 & 1.514 \\ 
$G_5$& CP & 0.979 & 0.977 & 0.985 & 0.959 & 0.961 & 0.967 \\  
&               & (0.933) & (0.930) & (0.935) & (0.920) & (0.915) & (0.910) \\ 
& $h^*$ & 0.772 & 0.752 & 1.088 & 0.885 & 1.112 & 1.648 \\ 
\hline
\end{tabular}
}
\caption{{Coverage probabilities (CP) for nominal $95\%$ confidence regions of the difference between two small area means for $k=2$ and $\D_i$-pattern (a).
(the corrected method in the first line and the naive method in parentheses)}}
\label{tab:cida}
\end{center}
\end{table}

\section{Appendix: proofs}

\subsection{Proof of Lemma \ref{lem:1}} 
The covariance of $\P \y$ and $\bbeh(\bPsi)$ is
\begin{align*}
E[&\P \y(\bbeh(\bPsi)-\bbe)^\top] \{\X^\top (\I_m\otimes \bPsi+\D)^{-1}\X\} \\
&=E[(\y-\X\bbeh^{OLS})(\y-\X\bbe)^\top](\I_m\otimes \bPsi+\D)^{-1}\X
\\
&=\Big[(\I_m\otimes \bPsi+\D) - \X \{\X^\top (\I_m\otimes \bPsi+\D)^{-1}\X\}^{-1}\X^\top \Big](\I_m\otimes \bPsi+\D)^{-1}\X
\\
&=\zero.
\end{align*}
This implies that $\bbeh(\bPsi)$ is independent of $\P \y$.
Next, we prove that $\bthh_a^{EB}-\bthh_a(\bPsi)$ is a function of $\P \y$.
From (H2), $\bPsih$ is a function of $\P \y$.
Rewrite $\bthh_a^{EB}$ as $\bthh_a^{EB}(\bPsih(\y),\y)$, $\bthh_a(\bPsi)$ as $\bthh_a(\bPsi,\y)$ and $\bbeh(\bPsi)$ as  $\bbeh(\bPsi,\y)$.
Since $\bbeh(\bPsi,\y+\X\T)=\bbeh(\bPsi,\y)+\T$, from (\ref{eqn:Bayes}) and (\ref{eqn:EB}), we have
\begin{align*}
&\bthh_a^{EB}(\bPsih(\y+\X\T),\y+\X\T)-\bthh_a(\bPsi(\y+\X\T),\y+\X\T)\\
=&\bthh_a^{EB}(\bPsih(\y),\y+\X\T)-\bthh_a(\bPsi(\y),\y+\X\T)\\
=&\y_a +\X_a\T - \D_a(\bPsih(\y)+\D_a)^{-1}\big\{\y_a+\X_a\T-\X_a\bbeh(\bPsih(\y))-\X_a\T\big\}\\
&-\y_a -\X_a\T + \D_a(\bPsi+\D_a)^{-1}\big\{\y_a+\X_a\T-\X_a\bbeh(\bPsi)-\X_a\T\big\}\\
=&\y_a  - \D_a(\bPsih(\y)+\D_a)^{-1}\big\{\y_a-\X_a\bbeh(\bPsih(\y))\big\}-\y_a  + \D_a(\bPsi+\D_a)^{-1}\big\{\y_a-\X_a\bbeh(\bPsi)\big\}\\
=&
\bthh_a^{EB}(\bPsih(\y),\y)-\bthh_a(\bPsi(\y),\y).
\end{align*}
Thus, $\bthh_a^{EB}-\bthh_a(\bPsi)$ is translation invariant, which implies that $\bthh_a^{EB}-\bthh_a(\bPsi)$ is a function of $\P \y$.
Hence, $\bthh_a^{EB}-\bthh_a(\bPsi)$ is independent of $\bbeh(\bPsi)$.
\hfill$\Box$

\subsection{Proof of Lemma \ref{lem:msemest}}

For the proof of part (2), note that 
\begin{equation}
(\bPsih+\D_i)^{-1}=(\bPsi+\D_i)^{-1}-(\bPsi+\D_i)^{-1}(\bPsih-\bPsi)(\bPsih+\D_i)^{-1}.
\label{eqn:identity}
\end{equation}
Then, $G_{1a}(\bPsih)$ is rewritten as 
\begin{align}
\G_{1a}(\bPsih)=& (\bPsih^{-1}+\D_a^{-1})^{-1}
=\D_a - \D_a(\bPsih+\D_a)^{-1}\D_a
\non\\
=&
\G_{1a}(\bPsi) + \D_a(\bPsi+\D_a)^{-1}(\bPsih-\bPsi)(\bPsi+\D_a)^{-1}\D_a
\label{eqn:BC0}\\
&- \D_a(\bPsi+\D_a)^{-1}(\bPsih-\bPsi)(\bPsi+\D_a)^{-1}(\bPsih-\bPsi)(\bPsi+\D_a)^{-1}\D_a+O_p(m^{-3/2}),\non
\end{align}
which implies that $E[\G_{1a}(\bPsih)]=\G_{1a}(\bPsi) - \G_{3a}(\bPsi)+O(m^{-3/2})$.

\medskip
For the proof of part (1), it is noted that
\begin{align*}
\bthh_a^{EB}-\bthh_a(\bPsi)
=& \D_a\{(\bPsi+\D_a)^{-1}-(\bPsih+\D_a)^{-1}\}(\y_a-\X_a\bbe)
+\D_a(\bPsih+\D_a)^{-1}\X_a\{\bbeh(\bPsih)-\bbe\}
\\
&-\D_a(\bPsi+\D_a)^{-1}\X_a\{\bbeh(\bPsi)-\bbe\}.
\end{align*}
Using the equation in (\ref{eqn:identity}),
we can see that
\begin{align*}
\D_a&\{(\bPsi+\D_a)^{-1}-(\bPsih+\D_a)^{-1}\}(\y_a-\X_a\bbe)\\
=& \D_a(\bPsi+\D_a)^{-1}(\bPsih-\bPsi)(\bPsih+\D_a)^{-1}(\y_a-\X_a\bbe)
\\
=&\D_a(\bPsi+\D_a)^{-1}(\bPsih-\bPsi)(\bPsi+\D_a)^{-1}(\y_a-\X_a\bbe)+O_p(m^{-1})
\end{align*}
and
\begin{align*}
\D_a&(\bPsih+\D_a)^{-1}\X_a\{\bbeh(\bPsih)-\bbe\}\\
=&\D_a(\bPsi+\D_a)^{-1}\X_a\{\bbeh(\bPsih)-\bbe\}
- \D_a(\bPsi+\D_a)^{-1}(\bPsih-\bPsi)(\bPsih+\D_a)^{-1}\X_a\{\bbeh(\bPsih)-\bbe\}
\\
=&\D_a(\bPsi+\D_a)^{-1}\X_a\{\bbeh(\bPsih)-\bbe\}+O_p(m^{-1}).
\end{align*}
Thus, we have
\begin{align*}
\bthh_a^{EB}-\bthh_a(\bPsi)
=&\D_a(\bPsi+\D_a)^{-1}(\bPsih-\bPsi)(\bPsi+\D_a)^{-1}(\y_a-\X_a\bbe)\\
&
+ \D_a(\bPsi+\D_a)^{-1}\X_a\{\bbeh(\bPsih)-\bbeh(\bPsi)\}+O_p(m^{-1})\\
=&I_1+I_2 + O_p(m^{-1}). \quad \text{(say)}
\end{align*}
For $I_2$, it is noted that
\begin{align*}
\bbeh&(\bPsih)-\bbeh(\bPsi)\\
=& \Big[\Big\{\sum_{j=1}^m \X_j^\top(\bPsih+\D_j)^{-1}\X_j\Big\}^{-1}-\Big\{\sum_{j=1}^m \X_j^\top(\bPsi+\D_j)^{-1}\X_j\Big\}^{-1}\Big] 
\\
&\times \sum_{i=1}^m \X_i^\top(\bPsih+\D_i)^{-1}(\y_i-\X_i\bbe)
\\
&+\Big\{\sum_{j=1}^m \X_j^\top(\bPsi+\D_j)^{-1}\X_j\Big\}^{-1}\sum_{i=1}^m \X_i^\top\Big\{(\bPsih+\D_i)^{-1}-(\bPsi+\D_i)^{-1}\Big\}(\y_i-\X_i\bbe)\\
=&I_{21}+I_{22}.
\end{align*}
We can evaluate $I_{21}$ as
\begin{align*}
I_{21}=&\Big\{\sum_{j=1}^m \X_j^\top(\bPsi+\D_j)^{-1}\X_j\Big\}^{-1}\Big\{\sum_{i=1}^m \X_i^\top (\bPsi+\D_i)^{-1}(\bPsih-\bPsi)(\bPsi+\D_i)^{-1}\X_i \Big\} \{\bbeh(\bPsih)-\bbe\}
\\
=& O_p(m^{-1}),
\end{align*}
because $\sum_{j=1}^m \X_j^\top(\bPsi+\D_j)^{-1}\X_j=O(m)$, $\sum_{i=1}^m \X_i^\top (\bPsi+\D_i)^{-1}(\bPsih-\bPsi)(\bPsi+\D_i)^{-1}\X_i =O_p(m^{1/2})$ and $\bbeh(\bPsih)-\bbe=O_p(m^{-1/2})$.
We next estimate $I_{22}$ as
\begin{align*}
I_{22}=&
- \Big\{\sum_{j=1}^m \X_j^\top(\bPsi+\D_j)^{-1}\X_j\Big\}^{-1}
\Big\{\sum_{i=1}^m \X_i^\top\A(\bPsih,\D_i)\X_i\Big\}\\
&\times \Big\{\sum_{i=1}^m \X_i^\top \A(\bPsih,\D_i)\X_i\Big\}^{-1}\sum_{i=1}^m \X_i^\top\A(\bPsih,\D_i)(\y_i-\X_i\bbe)
\end{align*}
for $\A(\bPsih,\D_i)=(\bPsih+\D_i)^{-1}(\bPsih-\bPsi)(\bPsi+\D_i)^{-1}$.
It can be seen that $I_{22}=O_p(m^{-1})$ from the same arguments as in $I_{21}$.
Thus, it follows that $I_2=O_p(m^{-1})$.
Hence, we have
\begin{align*}
E[\{&\bthh_a^{EB}-\bthh_a(\bPsi)\}\{\bthh_a^{EB}-\bthh_a(\bPsi)\}^\top]
\\
=&\D_a(\bPsi+\D_a)^{-1}E\Big[(\bPsih-\bPsi)(\bPsi+\D_a)^{-1}(\y_a-\X_a\bbe)
(\y_a-\X_a\bbe)^\top (\bPsi+\D_a)^{-1}(\bPsih-\bPsi)\Big]
\\
&\times  (\bPsi+\D_a)^{-1}\D_a + O(m^{-3/2}).
\end{align*}
Let $\bPsih_{(-a)}$ be an estimator of $\bPsi$ from the data except the $a$th area.
If we add or remove the data of one area in the estimation of $\bPsi$, there is a negligible change in the value of the above expectation since $\bPsih-\bPsih_{(-a)}=O_p(m^{-1})$.
Thus, we have
\begin{align*}
E[\{&\bthh_a^{EB}-\bthh_a(\bPsi)\}\{\bthh_a^{EB}-\bthh_a(\bPsi)\}^\top]
\\
=&\D_a(\bPsi+\D_a)^{-1}E\Big[(\bPsih_{(-a)}-\bPsi)(\bPsi+\D_a)^{-1}(\y_a-\X_a\bbe)
(\y_a-\X_a\bbe)^\top (\bPsi+\D_a)^{-1}(\bPsih_{(-a)}-\bPsi)\Big]
\\
&\times  (\bPsi+\D_a)^{-1}\D_a + O(m^{-3/2})
\\
=&\D_a(\bPsi+\D_a)^{-1}E\Big[(\bPsih_{(-a)}-\bPsi)(\bPsi+\D_a)^{-1}(\bPsih_{(-a)}-\bPsi)\Big](\bPsi+\D_a)^{-1}\D_a + O(m^{-3/2})
\\
=&
\D_a(\bPsi+\D_a)^{-1}E\Big[(\bPsih-\bPsi)(\bPsi+\D_a)^{-1}(\bPsih-\bPsi)\Big](\bPsi+\D_a)^{-1}\D_a + O(m^{-3/2}),
\end{align*}
which is equal to $\G_{1a}(\bPsi)+O(m^{-3/2})$, where the second equation follows from the independence of the data of different areas, and the the third equation follows form the same reason mentioned above.
\hfill$\Box$

\subsection{Proof of Theorem \ref{thm:cr}}
Let $\z_a=\H_a^{-1/2}(\bPsi)(\bthh_a^{EB}-\bth_a-\bthh_a^{EB}+\bthh_a(\bPsi))$.
From Lemma \ref{lem:cond}, the conditional distribution of $\z_a$ given $\P  \y$ is $\z_a \sim \Nc_k (\zero, \I_k)$, and the mahalanobis distance is approximated as
\begin{align}
&(\bthh_a^{EB}-\bth_a)^\top\H_a^{-1}(\bPsih)(\bthh_a^{EB}-\bth_a)
\non\\
=&
\z_a^\top\H_a^{1/2}(\bPsi)\H_a^{-1}(\bPsih)\H_a^{1/2}(\bPsi)\z_a
+2(\bthh_a^{EB}-\bthh_a(\bPsi))^\top\H_a^{-1}(\bPsih)\H_a^{1/2}(\bPsi)\z_a
\non\\
&+(\bthh_a^{EB}-\bthh_a(\bPsi))^\top\H_a^{-1}(\bPsih)(\bthh_a^{EB}-\bthh_a(\bPsi))
\non\\
=&
\z_a^\top\Big[\I_k-\H_a^{-1/2}(\bPsi)(\H_a(\bPsih)-\H_a(\bPsi))\H_a^{-1/2}(\bPsi)
\non\\
&\qquad+\H_a^{-1/2}(\bPsi)(\G_{1a}(\bPsih)-\G_{1a}(\bPsi))\H_a^{-2}(\bPsi)(\G_{1a}(\bPsih)-\G_{1a}(\bPsi))\H_a^{-1/2}(\bPsi)\Big]\z_a
\non\\
&+2(\bthh_a^{EB}-\bthh_a(\bPsi))^\top\H_a^{-1}(\bPsih)\H_a^{1/2}(\bPsi)\z_a\non\\
&
+(\bthh_a^{EB}-\bthh_a(\bPsi))^\top\H_a^{-1}(\bPsih)(\bthh_a^{EB}-\bthh_a(\bPsi))
+o(m^{-1})
\non\\
=&
\z_a^\top(\I_k-\G_{12a}(\bPsih))\z_a+2\g_{2a}(\bPsih)^\top\z_a+g_{3a}(\bPsih)+o(m^{-1}),
\label{eqn:mdap}
\end{align}
where 
\begin{align*}
\G_{12a}(\bPsih)=&
\H_a^{-1/2}(\bPsi)(\H_a(\bPsih)-\H_a(\bPsi))\H_a^{-1/2}(\bPsi)\\
&
-\H_a^{-1/2}(\bPsi)(\G_{1a}(\bPsih)-\G_{1a}(\bPsi))\H_a^{-2}(\bPsi)(\G_{1a}(\bPsih)-\G_{1a}(\bPsi))\H_a^{-1/2}(\bPsi),
\\
\g_{2a}(\bPsih)^\top=&(\bthh_a^{EB}-\bthh_a(\bPsi))^\top\H_a^{-1}(\bPsih)\H_a^{1/2}(\bPsi),\\
g_{3a}(\bPsih)=&
(\bthh_a^{EB}-\bthh_a(\bPsi))^\top\H_a^{-1}(\bPsih)(\bthh_a^{EB}-\bthh_a(\bPsi)).
\end{align*}
From (\ref{eqn:mdap}), the characteristic function $\fai(t)=E[\exp\{it (\bthh_a^{EB}-\bth_a)^\top\H_a^{-1}(\bPsih)(\bthh_a^{EB}-\bth_a)\}]$ is approximated as 
\begin{align*}
\fai(t)=&E\exp\Big(it\{\z_a^\top(\I_k-\G_{12a}(\bPsih))\z_a+2\g_{2a}(\bPsih)^\top\z_a+g_{3a}(\bPsih)\}\Big)+o(m^{-1})
\\
=&
E\Big[e^{it\z_a^\top\z_a} \Big\{ 1+it \{-\z_a^\top\G_{12a}(\bPsih)\z_a+2\g_{2a}(\bPsih)^\top\z_a+g_{3a}(\bPsih)\}
\\
&-{t^2 \over 2}  \{-\z_a^\top\G_{12a}(\bPsih)\z_a+2\g_{2a}(\bPsih)^\top\z_a+g_{3a}(\bPsih)\}^2\Big\} \Big]+o(m^{-1})
\\
=&
E\Big[e^{it\z_a^\top\z_a} \Big\{ 1+it \{-\z_a^\top\G_{12a}(\bPsih)\z_a+2\g_{2a}(\bPsih)^\top\z_a+g_{3a}(\bPsih)\}
\\
&-{t^2 \over 2}  \{(\z_a^\top\G_{12a}(\bPsih)\z_a)^2+4\z_a^\top\g_{2a}(\bPsih)\g_{2a}(\bPsih)^\top\z_a-4\z_a^\top\G_{12a}(\bPsih)\z_a\g_{2a}(\bPsih)^\top\z_a\}
\Big\} \Big]+o(m^{-1}),
\end{align*}
because $\G_{12a}(\bPsih)=O_p(m^{-1/2})$, $\g_{2a}(\bPsih)=O_p(m^{-1/2})$ and $g_{3a}(\bPsih)=O_p(m^{-1})$.
From the law of iterated expectations and the conditional normality of $\z_a$, the above equation reduces to
\begin{align*}
\fai(t)=&
E\Big[e^{it\z_a^\top\z_a} \Big\{ 1+it \{-\z_a^\top\G_{12a}(\bPsih)\z_a+g_{3a}(\bPsih)\}\\
&
\qquad\qquad-{t^2 \over 2}  \{(\z_a^\top\G_{12a}(\bPsih)\z_a)^2
+4\z_a^\top\g_{2a}(\bPsih)\g_{2a}(\bPsih)^\top\z_a\}
\Big\} \Big]
+o(m^{-1}).
\end{align*}

For some deterministic matrix $\A$ and $\z \sim \Nc_k (\zero, \I_k)$, it holds that
\begin{align*}
E\Big[e^{it\z^\top\z} \z^\top \A \z \Big]=&
(2\pi)^{-k/2}\int e^{-{(1-2it)\z^\top\z \over 2}}\z^\top \A \z d\z
=
(1-2it)^{-k/2-1} \tr (\A),
\\
E\Big[e^{it\z^\top\z} (\z^\top \A \z)^2 \Big]=&
(2\pi)^{-k/2}\int e^{-{(1-2it)\z^\top\z \over 2}}(\z^\top \A \z)^2 d\z
=
(1-2it)^{-k/2-2} (\tr ^2(\A)+2\tr (\A^2)).
\end{align*}
Using these equalities, from the law of iterated expectations, we have
\begin{align*}
\fai(t)=&
(1-2it)^{-k/2} \Big[ 1+it \Big\{-(1-2it)^{-1}\tr(E[\G_{12a}(\bPsih)])+E[g_{3a}(\bPsih)]\Big\}\\
&\qquad+{(it)^2 \over 2}  \Big\{(1-2it)^{-2} \{E[\tr^2( \G_{12a}(\bPsih))]+2\tr( E[\G_{12a}^2(\bPsih)]) \}
\\
&\qquad\qquad+(1-2it)^{-1}4\tr(E[\g_{2a}(\bPsih)\g_{2a}(\bPsih)^\top])\Big\}
\Big]+o(m^{-1}).
\end{align*}
For notational simplicity, let $C=E[\tr^2( \G_{12a}(\bPsih))]+2\tr( E[\G_{12a}^2(\bPsih)])$.
Let $s= (1-2it)^{-1}$, or $it=(s-1)/2s$.
Then,  $(1-2it)^{-k/2}\fai(t)-1$ can be written as
\begin{align}
&it \Big\{-(1-2it)^{-1}\tr(E[\G_{12a}(\bPsih)])+E[g_{3a}(\bPsih)]\Big\}\non\\
&+{(it)^2 \over 2}  \Big\{(1-2it)^{-2} C 
+(1-2it)^{-1}4E[\g_{2a}(\bPsih)^\top\g_{2a}(\bPsih)]\Big\}
\non\\
=&
{1\over 2s}\Big\{E[\g_{2a}(\bPsih)^\top\g_{2a}(\bPsih)]-E[g_{3a}(\bPsih)]\Big\}
\non\\
&+ \Big\{{1 \over 2}\tr( E[\G_{12a}(\bPsih)])+{1 \over 2}E[g_{3a}(\bPsih)] +{C \over 8}- E[\g_{2a}(\bPsih)^\top\g_{2a}(\bPsih)]\Big\}
\non\\
&+\Big\{ -{1 \over 2}\tr( E[\G_{12a}(\bPsih)])-{C \over 4}
+{1 \over 2}E[\g_{2a}(\bPsih)^\top\g_{2a}(\bPsih)]\Big\}s
+{C \over 8} s^2
+o(m^{-1}).
\label{eqn:pol}
\end{align}
which is a second-order polynomial of $s$. 

\medskip
We shall evaluate the moments in (\ref{eqn:pol}).
First, $\G_{12a}(\bPsih)$ can be expanded as
\begin{align}
\G_{12a}(\bPsih)
=&
\H_a^{-1/2}(\bPsi)(\H_a(\bPsih)-\H_a(\bPsi))\H_a^{-1/2}(\bPsi)
\label{eqn:G12ap}
\\
&-\H_a^{-1/2}(\bPsi)(\G_{1a}(\bPsih)-\G_{1a}(\bPsi))\H_a^{-2}(\bPsi)(\G_{1a}(\bPsih)-\G_{1a}(\bPsi))\H_a^{-1/2}(\bPsi)
+o_p(m^{-1}).
\non
\end{align}
From Lemma \ref{lem:msemest}, the expectation of the first term in (\ref{eqn:G12ap}) is $-\H_a^{-1/2}(\bPsi)\G_{3a}(\bPsi)\H_a^{-1/2}(\bPsi)+o(m^{-1})$, so that 
$$
E[\G_{12a}(\bPsih)]=-\H_a^{-1/2}(\bPsi)\G_{3a}(\bPsi)\H_a^{-1/2}(\bPsi)
- E[\K_a(\bPsih)\H_a^{-1}(\bPsi)\K_a(\bPsih)] + o(m^{-1}),
$$
for $\K_a(\bPsih)=\H_a^{-1/2}(\bPsi)(\G_{1a}(\bPsih)-\G_{1a}(\bPsi))\H_a^{-1/2}$.
Thus,
\begin{equation}
\tr(E[\G_{12a}(\bPsih)])=- B_3+2B_1 + o(m^{-1}),
\label{eqn:pcr1}
\end{equation}
for $B_1$ and $B_3$ defined in (\ref{eqn:B12}).
Noting that the first term in (\ref{eqn:G12ap}) is of order $O(m^{-1/2})$ and the second term is of order $O(m^{-1})$, we can expand $\G_{12a}^2(\bPsih)$ and $\tr^2(\G_{12a}(\bPsih))
$ as
\begin{align*}
\G_{12a}^2(\bPsih)
=&
\H_a^{-1/2}(\bPsi)(\G_{1a}(\bPsih)-\G_{1a}(\bPsi))\H_a^{-1}(\bPsi)(\G_{1a}(\bPsih)-\G_{1a}(\bPsi))\H_a^{-1/2}(\bPsi)+o_p(m^{-1}),
\\
\tr^2(\G_{12a}(\bPsih))
=&
\tr^2(\H_a^{-1/2}(\bPsi)(\G_{1a}(\bPsih)-\G_{1a}(\bPsi))\H_a^{-1/2}(\bPsi))+o_p(m^{-1}),
\end{align*}
which lead to $E[\G_{12a}^2(\bPsih)]=E[ \{ \K_a(\bPsih) \}^2] +o(m^{-1})$ and $E[\tr^2(\G_{12a}(\bPsih))]=E[\tr^2(\K_a(\bPsih))]+o(m^{-1})$.
Thus,
\begin{equation}
C=E[\tr^2(\K_a(\bPsih))]+2E[ \{ \K_a(\bPsih) \}^2]+o(m^{-1})
=-8B_2 +o(m^{-1}),
\label{eqn:pcr2}
\end{equation}
for $B_2$ defined in (\ref{eqn:B12}).
It can be also observed that
\begin{align*}
\g_{2a}(\bPsih)^\top\g_{2a}(\bPsih)
=&
(\bthh_a^{EB}-\bthh_a(\bPsi))^\top\H_a^{-1}(\bPsih)\H_a(\bPsi)\H_a^{-1}(\bPsih)(\bthh_a^{EB}-\bthh_a(\bPsi))
\\
=&
(\bthh_a^{EB}-\bthh_a(\bPsi))^\top\H_a^{-1}(\bPsi)(\bthh_a^{EB}-\bthh_a(\bPsi))+o_p(m^{-1}),
\\
g_{3a}(\bPsih)
=&
(\bthh_a^{EB}-\bthh_a(\bPsi))^\top\H_a^{-1}(\bPsih)(\bthh_a^{EB}-\bthh_a(\bPsi))
\\
=&
(\bthh_a^{EB}-\bthh_a(\bPsi))^\top\H_a^{-1}(\bPsi)(\bthh_a^{EB}-\bthh_a(\bPsi))+o_p(m^{-1}),
\end{align*}
both of which lead to 
\begin{equation}
E[\g_{2a}(\bPsih)^\top\g_{2a}(\bPsih)]=E[ g_{3a}(\bPsih)]
=
\tr(\H_a^{-1}(\bPsi)\G_{3a}(\bPsi))+o(m^{-1})=B_3+o(m^{-1}).
\label{eqn:pcr3}
\end{equation}

Combining (\ref{eqn:pcr1}), (\ref{eqn:pcr2}) and (\ref{eqn:pcr3}), we can see that the constant term and the coefficient of $s^2$ in (\ref{eqn:pol})  are $B_1-B_3-B_2$ and $-B_2$ given in (\ref{eqn:B12}), respectively.
Thus, the characteristic function of $(\bthh_a^{EB}-\bth_a)^\top\H_a^{-1}(\bPsih)(\bthh_a^{EB}-\bth_a)$ can be written as
\begin{align*}
\fai(t)=(1-2it)^{-k/2}(1+B_1-B_3-B_2+(-B_1+B_3+2B_2)s-B_2s^2)+o(m^{-1}).
\end{align*}
From the fact that the characteristic function of the chi-squared distribution with degrees of freedom $k+2h$ is given by $(1-2it)^{-k/2-h}=(1-2it)^{-k/2}s^h$, it follows that the asymptotic expansion of the cumulative distribution function of $(\bthh_a^{EB}-\bth_a)^\top\H_a^{-1}(\bPsih)(\bthh_a^{EB}-\bth_a)$ is 
\begin{align*}
F_{k}(x)+(B_1-B_3-B_2)F_{k}(x)+(-B_1+B_3+2B_2)F_{k+2}(x)-B_2F_{k+4}(x)+o(m^{-1}),
\end{align*}
where $F_{k}(x)$ is the cumulative distribution function of the chi-squared distribution with degrees of freedom $k$.
Note that $F_{k+r-2}(x)-F_{k+r}(x)=2f_{k+r}(x)$, where $f_k(x)$ is the density function of the chi-squared distribution with degrees of freedom $k$.
Then, it is expressed as
\begin{align*}
P((\bthh_a^{EB}&-\bth_a)^\top\H_a^{-1}(\bPsih)(\bthh_a^{EB}-\bth_a)\leq x)\\
=&F_{k}(x)+2(B_1-B_3-B_2)f_{k+2}(x)+2B_2f_{k+4}(x)+o(m^{-1}),
\end{align*}
which proves Theorem \ref{thm:cr}.
\hfill$\Box$

\subsection{Proof of Lemma \ref{lem:MSE}}
From Proposition \ref{prp:pd}, it is sufficient to show this approximation for $\bPsih^{PR}$ instead of $\bPsih_{(A)}^{PR}$.
It is noted that $\bPsih^{PR}-\bPsi$ is approximated as
\begin{align}
\bPsih^{PR}-\bPsi=&
{1\over m}\sum_{i=1}^m\{(\y_i-\X_i\bbe)(\y_i-\X_i\bbe)^\top - (\bPsi+\D_i)\}
+O_p(m^{-1}),
\label{eqn:Psiap}
\end{align}
which is used to evaluate 
\begin{align*}
E\Big[&(\bPsih-\bPsi)(\bPsi+\D_a)^{-1}(\y_a-\X_a\bbe)
(\y_a-\X_a\bbe)^\top (\bPsi+\D_a)^{-1}(\bPsih-\bPsi)\Big]
\\
=&
{1\over m^2}\sum_{i=1}^m\sum_{j=1}^m 
E\Big[ \{\u_i\u_i^\top-(\bPsi+\D_i)\}(\bPsi+\D_a)^{-1}\u_a\u_a^\top (\bPsi+\D_a)^{-1}
\{\u_j\u_j^\top-(\bPsi+\D_i)\}\Big]\\
&+O(m^{-3/2})\non\\
=&
{1\over m^2}\sum_{i=1}^m
E\Big[ \{\u_i\u_i^\top-(\bPsi+\D_i)\}(\bPsi+\D_a)^{-1}\u_a\u_a^\top (\bPsi+\D_a)^{-1}
\{\u_i\u_i^\top-(\bPsi+\D_i)\}\Big]\\
&+O(m^{-3/2}),
\end{align*}
because $E\Big[ \{\u_i\u_i^\top-(\bPsi+\D_i)\}(\bPsi+\D_a)^{-1}\u_a\u_a^\top (\bPsi+\D_a)^{-1}\{\u_j\u_j^\top-(\bPsi+\D_i)\}\Big]=\zero$ for $i\not= j$.
Letting $\z_i=(\bPsi+\D_i)^{-1/2}\u_i$, we can see that  $\z_i\sim \Nc_k(\zero, \I_k)$.
Then,
\begin{align*}
{1\over m^2}&\sum_{i=1}^m
E\Big[ \{\u_i\u_i^\top-(\bPsi+\D_i)\}(\bPsi+\D_a)^{-1}\u_a\u_a^\top (\bPsi+\D_a)^{-1}
\{\u_i\u_i^\top-(\bPsi+\D_i)\}\Big]\\
=&
{1\over m^2}\sum_{i\not= a}
(\bPsi+\D_i)^{1/2}E\Big[ (\z_i\z_i^\top-\I)\B\z_a\z_a^\top \B^\top (\z_i\z_i^\top-\I)\Big](\bPsi+\D_i)^{1/2} + O(m^{-2}),
\end{align*}
for $\B=(\bPsi+\D_i)^{1/2}(\bPsi+\D_a)^{-1/2}$.
Let $\C=\B\B^\top=(\bPsi+\D_i)^{1/2} (\bPsi+\D_a)^{-1} (\bPsi+\D_i)^{1/2}$.
For $i\not= a$, 
\begin{align*}
E[& (\z_i\z_i^\top-\I)\B\z_a\z_a^\top \B^\top (\z_i\z_i^\top-\I)]
\\
=&E[ \z_i\z_i^\top\B\z_a\z_a^\top \B^\top \z_i\z_i^\top+\B\z_a\z_a^\top \B^\top
-\z_i\z_i^\top\B\z_a\z_a^\top \B^\top -\B\z_a\z_a^\top \B^\top \z_i\z_i^\top]
\\
=&
E[ \z_i\z_i^\top\C \z_i\z_i^\top-\C]
=\C+(\tr\C)\I_k,
\end{align*}
because $E[\z_i\z_i^\top\C \z_i\z_i^\top]=2\C+(\tr\C)\I_k$.
Thus,
\begin{align*}
{1\over m^2}&\sum_{i\not= a}
(\bPsi+\D_i)^{1/2}\{ \C+(\tr\C)\I_k\}(\bPsi+\D_i)^{1/2}
\\
=&
{1\over m^2}\sum_{i=1}^m
(\bPsi+\D_i)^{1/2}\{ \C+(\tr\C)\I_k\}(\bPsi+\D_i)^{1/2} + O(m^{-2}),
\end{align*}
which leads to the expression in (\ref{eqn:G3}) from Lemma \ref{lem:msemest} (1).
\hfill$\Box$

\subsection{Proof of Lemma \ref{lem:moment}}

From Proposition \ref{prp:pd}, it is sufficient to show this approximation for $\bPsih^{PR}$ instead of $\bPsih_{(A)}^{PR}$.
For some deteministic matrix $\A$ and multivariate standard normal variables $\z_i$, $i=1,\ldots,k$, $E[(\bPsih^{PR}-\bPsi)\A(\bPsih^{PR}-\bPsi)]$ is, from (\ref{eqn:Psiap}),  approximated as
\begin{align*}
E[(&\bPsih^{PR}-\bPsi)\A(\bPsih^{PR}-\bPsi)]\\
=&{1\over m^2}\sum_{i=1}^m\sum_{j=1}^m E\Big[\Big\{(\y_i-\X_i\bbe)(\y_i-\X_i\bbe)^\top - (\bPsi+\D_i)\Big\}\A
\\
&\times \Big\{(\y_j-\X_j\bbe)(\y_j-\X_j\bbe)^\top - (\bPsi+\D_j)\Big\}\Big] + O(m^{-3/2})
\\
=&
{1\over m^2}\sum_{i=1}^m\sum_{j=1}^m (\bPsi+\D_i)^{1/2}
E[( \z_i\z_i^\top-\I) \C_i ( \z_j\z_j^\top-\I)]  (\bPsi+\D_j)^{1/2} + O(m^{-3/2}),
\end{align*}
for $\C_i=(\bPsi+\D_i)^{1/2} \A (\bPsi+\D_i)^{1/2}$.
For $i\not= j$, $E[( \z_i\z_i^\top-\I) \C_i ( \z_j\z_j^\top-\I)] =0$, we have 
$$
\sum_{i=1}^m\sum_{j=1}^m E[( \z_i\z_i^\top-\I) \C_i ( \z_j\z_j^\top-\I)] 
=\sum_{i=1}^m E[\z_i\z_i^\top\C_i\z_i\z_i^\top-\C_i].
$$
Because $E[\z_i\z_i^\top\C_i\z_i\z_i^\top]=2\C_i+(\tr\C_i)\I_k$, it is concluded that 
\begin{align*}
E[(\bPsih^{PR}-\bPsi)\A(\bPsih^{PR}-\bPsi)]={1\over m^2}\sum_{i=1}^m((\bPsi+\D_i) \A (\bPsi+\D_i)+\tr( \A (\bPsi+\D_i))(\bPsi+\D_i)).
\end{align*}

Using this equality, we have 
\begin{align*}
&E[(\bPsih^{PR}-\bPsi)(\bPsi+\D_a)^{-1}\D_a\H_a^{-2}(\bPsi)\D_a(\bPsi+\D_a)^{-1}(\bPsih^{PR}-\bPsi)]
\\
=&
{1 \over m^2} \sum_{i=1}^m \Big\{ (\bPsi+\D_i)(\bPsi+\D_a)^{-1}\D_a\H_a^{-2}(\bPsi)\D_a(\bPsi+\D_a)^{-1} (\bPsi+\D_i)
\\
&+ \tr((\bPsi+\D_a)^{-1}\D_a\H_a^{-2}(\bPsi)\D_a(\bPsi+\D_a)^{-1} (\bPsi+\D_i))(\bPsi+\D_i)\Big\}
+o(m^{-1}),
\end{align*}
and
\begin{align*}
&E[(\G_{1a}(\bPsih^{PR})-\G_{1a}(\bPsi))\H_a^{-1}(\bPsi)(\G_{1a}(\bPsih^{PR})-\G_{1a}(\bPsi))]
\\
=&
\D_a\H_a^{-1}(\bPsi) \Big[{1 \over m^2} \sum_{i=1}^m \Big\{ (\bPsi+\D_i)(\bPsi+\D_a)^{-1}\D_a\H_a^{-1}(\bPsi)\D_a(\bPsi+\D_a)^{-1} (\bPsi+\D_i)
\\
&+ \tr((\bPsi+\D_a)^{-1}\D_a\H_a^{-1}(\bPsi)\D_a(\bPsi+\D_a)^{-1} (\bPsi+\D_i))(\bPsi+\D_i)\Big\}\Big](\bPsi+\D_a)^{-1}\D_a
+o(m^{-1}),
\end{align*}
which leads to the first and third expressions in the lemma.
From (\ref{eqn:BC0}), 
\begin{align*}
E[\tr^2(\G_{12a}(\bPsih^{PR}))]
=&
E[\tr^2(\H_a^{-1}(\bPsi)(\G_{1a}(\bPsih^{PR})-\G_{1a}(\bPsi)))]
+o(m^{-1})
\\
=&
E[\tr^2((\bPsi+\D_a)^{-1}\D_a\H_a^{-1}(\bPsi)\D_a(\bPsi+\D_a)^{-1} (\bPsih^{PR}-\bPsi))]
+o(m^{-1}).
\end{align*}
Letting $\u_i=\y_i-\X_i\bbe$, we can see that
\begin{align*}
&\tr((\bPsi+\D_a)^{-1}\D_a\H_a^{-1}(\bPsi)\D_a(\bPsi+\D_a)^{-1} (\bPsih^{PR}-\bPsi))
\\
=&
{1 \over m}\sum_{i=1}^{m}\tr((\bPsi+\D_a)^{-1}\D_a\H_a^{-1}(\bPsi)\D_a(\bPsi+\D_a)^{-1} (\bPsi+\D_i)^{1/2}(\u_i\u_i^\top-\I_k)(\bPsi+\D_i)^{1/2})
\\
&+o(m^{-1})
\\
=&
{1 \over m}\sum_{i=1}^{m}\tr((\bPsi+\D_a)^{-1}\D_a\H_a^{-1}(\bPsi)\D_a(\bPsi+\D_a)^{-1} (\bPsi+\D_i)^{1/2}\u_i\u_i^\top(\bPsi+\D_i)^{1/2})
\\
&-{1 \over m}\sum_{i=1}^{m}\tr((\bPsi+\D_a)^{-1}\D_a\H_a^{-1}(\bPsi)\D_a(\bPsi+\D_a)^{-1} (\bPsi+\D_i))
+o(m^{-1}).
\end{align*}
Thus we have
\begin{align*}
E[\tr^2(\G_{12a}(\bPsih^{PR}))]
=
{2 \over m^2}\sum_{i=1}^{m}\tr\Big(((\bPsi+\D_a)^{-1}\D_a\H_a^{-1}(\bPsi)\D_a(\bPsi+\D_a)^{-1} (\bPsi+\D_i))^2\Big)
+o(m^{-1}),
\end{align*}
which leads to the expression in the lemma.
\hfill$\Box$

\section*{Acknowledgments}
Research of the second author was supported in part by Grant-in-Aid for Scientific Research  (15H01943 and 26330036) from Japan Society for the Promotion of Science.

\end{document}